\newif\ifwias\wiastrue
\newcount\hours\newcount\minutes
\def \SetTime{\hours=\time\global\divide\hours by 60\minutes=\hours
\multiply\minutes by 60\advance\minutes by-\time\global\multiply\minutes by-1}
\SetTime
\def \now{\number\hours:\ifnum\minutes<10 0\fi\number\minutes}

\documentclass[a4paper,12pt,twoside,draft]{article}
\usepackage[english,german]{babel}
\usepackage{ifpdf}
\ifwias
\usepackage{wiasa4,wiaspreprint}
\renewtextwidthandheight{15cm}{9in}
\else
\usepackage{vmargin}
\setpapersize[portrait]{A4}
\setmarginsrb{30mm}{20mm}{30mm}{40mm}{10mm}{5mm}{10mm}{15mm}
\fi
\usepackage{amsfonts,amsmath,amstext,amsbsy,amsopn,amsthm}
\usepackage{amssymb} 
\usepackage{pifont} 
\usepackage{verbatim}
\usepackage{ifthen}
\ifpdf
\usepackage[OT1]{fontenc}
\usepackage[hyperindex,colorlinks=true,citecolor=blue]{hyperref}
\usepackage[pdftex]{graphicx}
\else
\usepackage[hyperindex,colorlinks=false,citecolor=blue]{hyperref}
\usepackage{graphicx}
\fi
\usepackage{fancyheadings}
\DeclareMathOperator{\INTERIOR}{interior}
\DeclareMathOperator{\dom}{dom}

\DeclareMathOperator{\supp}{supp}  %support of a function
\DeclareMathOperator{\trace}{tr}   %trace of a function on a boundary
\DeclareMathOperator*{\essinf}{vrai min}

\newcommand{\abs}[1]{\lvert#1\rvert}

\newcommand{\norm}[1]{\lVert#1\rVert}

\newcommand{\lrnorm}[1]{\left\lVert#1\right\rVert}  

\newcommand{\dual}[2]{\langle{#1}\,|\,{#2}\rangle}

\newcommand{\lrdual}[2]{\left\langle{#1}\,|\,{#2}\right\rangle}

\newcommand{\upp}[1]{{#1}^\bullet}
\newcommand{\low}[1]{{#1}_\bullet}

\newcommand{\xoplus}[1]{\operatorname{#1}}

\newcommand{\ta}{t}
\newcommand{\tb}{\check{t}}

\newcommand{\constem}{E_M}
\newcommand{\constgm}{G_M}
\newcommand{\constqm}{Q_M}
\newcommand{\constzm}{Z_M}
\newcommand{\constym}{Y_M}

\newcommand{\TNULL}{T_0}
\newcommand{\TEINS}{T_1}

\newcommand{\nlpsolution}{\mathcal{L}}

\newcommand{\id}{\mathbb{I}}

\newcommand{\e}{\mathrm{e}}
\newcommand{\exponential}[1]{\e^{#1}}

\newcommand{\der}{\,\mathrm{d}}
\newcommand{\df}{\stackrel{\scriptscriptstyle\mathrm{def}}{=}}
\newcommand{\textfrac}[2]{{#1}/{#2}}

\newcommand{\with}{\,:\,}
\newcommand{\from}[1]{:#1}

\newcommand{\interior}[1]{\INTERIOR({#1})}

\newcommand{\dive}{\nabla\cdot}
\newcommand{\grad}{\nabla}
\newcommand{\DIVE}{\nabla\cdot} %%\diag{\dive}
\newcommand{\GRAD}{\nabla} %%\diag{\gard}

\newcommand{\Lindensity}{\mathcal{N}}
\newcommand{\Linpoisson}{\mathcal{P}_0}

\newcommand{\fulla}[1]{\check{a}_{#1}}

\newcommand{\matrixA}{\operatorname{A}}

\newcommand{\domainA}{\mathcal{D}}

\newcommand{\Dirichlet}{\mathrm{D}}
\newcommand{\widehatDirichlet}{\widehat{\Dirichlet}}

\newcommand{\Neumann}{\Gamma}
\newcommand{\widehatNeumann}{\widehat{\Neumann}}

\newcommand{\embedcontinuously}{\hookrightarrow}

\newcommand{\subomega}{\omega}
\newcommand{\restrict}{\mathstrut_\downarrow}
\newcommand{\extension}{\mathstrut^\uparrow}
\newcommand{\reaction}{r}
\newcommand{\chempot}{\chi}
\newcommand{\doping}{\tilde{d}}
\newcommand{\immotile}{d}

\newcommand{\Varphi}{\widetilde{\varphi}}
\newcommand{\Fullphi}{\widetilde{\phi}}

\newcommand{\Dbvarphi}{\varphi_{\widehatDirichlet}}
\newcommand{\Rbvarphi}{\varphi_{\widehatNeumann}}     
\newcommand{\rbvarphi}{\varphi_\bullet}   %%% function on the domain
\newcommand{\dbvarphi}{\varphi_\circ}     %%% function on the domain

\newcommand{\Bvphi}[1]{\phi_{\Dirichlet,#1}} 
\newcommand{\bvphi}[1]{{\Phi}_{#1}} 
\newcommand{\iniphi}{\Phi^0}  % initial value, function on the domain 

\newcommand{\bvpsi}{\Psi}
\newcommand{\wildbvpsi}{\breve{\psi}} %%% function on the domain 

\newcommand{\cf}{see}
\newcommand{\ie}{i.e.}
\newcommand{\eg}{for instance}

\swapnumbers

\newtheorem{thm}{Theorem}[section]
\newtheorem{lem}[thm]{Lemma}%[section]
\newtheorem{cor}[thm]{Corollary}%[section]
\newtheorem{prop}[thm]{Proposition}%[section]
\theoremstyle{definition}
\newtheorem{defn}[thm]{Definition}%[section]
\newtheorem{assu}[thm]{Assumption}%[section]
\newtheorem{rem}[thm]{Remark}%[section]

\theoremstyle{remark}
\newtheorem*{acknowledgement}{Acknowledgement}  

\numberwithin{equation}{section}

\newcommand{\assuref}[1]{Assumption~\ref{#1}}
\newcommand{\thmref}[1]{Theorem~\ref{#1}}
\newcommand{\corref}[1]{Corollary~\ref{#1}}
\newcommand{\remref}[1]{Remark~\ref{#1}}
\newcommand{\proref}[1]{Proposition~\ref{#1}}
\newcommand{\secref}[1]{\S\ref{#1}}
\newcommand{\lemref}[1]{Lemma~\ref{#1}}

\newcommand{\defnref}[1]{Definition~\ref{#1}}

\newcommand{\thmtitle}[1]{\emph{#1}}

\begin{document}

\title{Classical solutions of
  % two-dimensional 
  drift--diffusion equations
  \\
  for semiconductor devices: the 2d case
% \title{Classical solutions
%   %% in Lebesgue spaces 
%   of two-dimensional drift--diffusion equations
%   for semiconductor devices
  \\{\normalsize\emph{Dedicated to 
      Herbert Gajewski, 
      Konrad Gr\"oger and 
      Klaus Zacharias}}}

\ifwias

\author{%
Hans-Christoph Kaiser,
Hagen Neidhardt\nofnmark{}
%\footnote{supported by the DFG under grant no.~RE 1480/2}, 
and 
Joachim Rehberg\nofnmark{}\footnote{%
\href{http://www.wias-berlin.de}{%
Weierstrass Institute for
Applied Analysis and Stochastics}\\
Mohrenstr. 39\\
10117 Berlin\\
Germany \\
\\
\begin{tabular}{ll}
E-Mail:
        &\href{mailto:kaiser@wias-berlin.de}{kaiser@wias-berlin.de}\\
        &\href{mailto:neidhardt@wias-berlin.de}{neidhardt@wias-berlin.de}\\
        &\href{mailto:rehberg@wias-berlin.de}{rehberg@wias-berlin.de}
\end{tabular}
}
}

\nopreprint{1189}   % Preprint-Nummer 
\nopreyear{2006} % Jahr des Preprints
% 35K45 Initial value problems for parabolic systems
% 35K50 Boundary value problems for parabolic systems
% 35K55 Nonlinear PDE of parabolic type
% 35K57 Reaction-diffusion equations
% 78A35 Motion of charged particles
\selectlanguage{english}% hier nicht veraendern, wichtig fuer Datumsformat 

\else

\author{Hans-Christoph Kaiser, Hagen Neidhardt, and Joachim Rehberg}

\fi

%\date{  D R A F T --- \Now\ --- D R A F T  }    % Datum fixieren
\date{December 15th, 2006}    % Datum fixieren

%% amsart only! abstract befor maketitle
%\begin{abstract} \end{abstract}
\maketitle
%% article and other classes: abstract after maketitle 

\begin{abstract}
  We regard drift--diffusion equations for semiconductor devices in
  Lebesgue spaces. To that end we reformulate the (generalized) van
  Roosbroeck system as an evolution equation for the potentials to the
  driving forces of the currents of electrons and holes.  This
  evolution equation falls into a class of quasi-linear parabolic
  systems which allow unique, local in time solution in certain
  Lebesgue spaces. In particular, it turns out that the divergence of
  the electron and hole current is an integrable function. Hence,
  Gauss' theorem applies, and gives the foundation for space
  discretization of the equations by means of finite volume schemes.
  Moreover, the strong differentiability of the electron and hole
  density in time is constitutive for the implicit time discretization
  scheme. 
  Finite volume discretization of space, and implicit time discretization are  
  accepted custom in engineering and scientific computing.
  --- This investigation puts special emphasis on non-smooth spatial
  domains, mixed boundary conditions, and heterogeneous material
  compositions, as required in electronic device simulation.
\end{abstract}

\newpage

\section[Introduction]{Introduction}
\label{sec:intro}

In 1950 van~Roosbroeck \cite{roosbroeck50} established a system of
partial differential equations describing the motion of electrons and
holes in a semiconductor device due to drift and diffusion within a
self-consistent electrical field. In 1964 Gummel \cite{gummel64}
published the first report on the numerical solution of these
drift--diffusion equations for an operating semiconductor device. From
that time on van~Roosbroeck's system has been the backbone of many a
model in semiconductor device simulation. The first papers devoted to
the mathematical analysis of van~Roosbroeck's system appeared in the
early seventies of the previous century \cite{mock72,mock74}; for a
historical synopsis and further references see \cite{gajewski93}.  In
1986 Gajewski and Gr{\"o}ger proved the global existence and
uniqueness of weak solutions under realistic physical and geometrical
conditions \cite{gajewski:groeger86}. The key for proving these
results and also for establishing stable numerical solving procedures
is the existence of a Lyapunov function for the van Roosbroeck system.
This solution theory entails restricting conditions on the models for
the recombination of electron--hole pairs, see
\cite[2.2.3]{gajewski93}, \cite[Ch.~5]{gajewski:groeger89},
\cite[Ch.~6]{gajewski:groeger90}, \cite{skrypnik:02}, and
\cite{pub:938}. In this paper we relax the condition on the reaction
terms in the equations considerably, up to the point that some
external control to the generation or annihilation of electrons or
holes can be applied individually. In particular, this aims at
radiative recombination of electron-hole pairs in semiconductor
lasers, and at the generation of electron-hole pairs in optoelectronic
detectors. Notwithstanding this generalization, we continue to use the
name van~Roosbroeck system for the model equations.

Van~Roosbroeck's system consists of current--continuity equations ---
one for electrons, another one for holes --- which are coupled to a
Poisson equation for the electrostatic potential, and comprise
generative terms, first of all recombination of electron--hole pairs.
The current--continuity equations can be viewed as quasi-linear
parabolic equations.  However, the natural formulation of balance laws
is in integral form
\begin{equation}
  \label{eq:balancelaw}
  \frac{\partial}{\partial t} 
  \int_\subomega u_k \der{x} 
  = 
  \int_{\partial\subomega} 
  {\nu}\cdot{j_k} \der{\sigma_\omega}
  +
  \int_{\subomega} \reaction_k \der{x}
  .
\end{equation}
Here $u_2$ and $u_1$ is the density of electrons and holes,
respectively, $j_k$ is the corresponding flux, and $\reaction_k$ is a
reaction term.  $\subomega$ is any (suitable) sub-domain of the whole
domain under consideration, $\nu$ the outer unit normal to the
boundary $\partial\subomega$ of $\subomega$ and $\sigma_\omega$ the
arc measure on $\partial \omega$.  In the weak formulation of the
balance law the boundary integral of the normal component of the
current is expressed as the volume integral of the divergence of the
corresponding current.
% seefor instance \cite[should one give references for that at
%all?]{gajewski:groeger90}, \cite{huenlich:95},
%\cite{glitzky:huenlich:97}.
                                %
Very little is known about the question whether the weak solutions also
satisfy the original balance law equations \eqref{eq:balancelaw}.
Obviously, this depends on the applicability of Gauss' theorem. So,
the problem is about the divergence of the currents in weak solutions
being functions --- not only distributions. In particular, this comes
to bear in the numerical treatment of van Roosbroeck's system. The
choice for space discretization of drift--diffusion equations is the
finite volume method, see \cite{gklnr:detectors}, which rests on
the original balance law formulation \eqref{eq:balancelaw} of the
equations.

In this paper we solve this problem for the spatially two-dimensional
van Roosbroeck system by showing that it admits a classical solution
in a suitably chosen Lebesgue space---at least locally in time. Aiming
at the inclusion of rather general recombination and generation
processes for electron-hole pairs we cannot expect global existence
anymore, and we cannot rely on a Lyapunov function. Instead we apply
local methods for quasi-linear evolution equations. To that end, we
rewrite van Roosbroeck's system as an evolution equation for the
electrochemical potentials of electrons and holes, and apply a
recently obtained result on quasi-linear parabolic equations in
Lebesgue spaces, see \cite{pub:765}. This yields a classical solution
of van Roosbroeck system locally in time with currents the divergence
of which is Lebesgue integrable to some exponent greater than one.
The strong differentiability of the electron and hole density in time
is constitutive for the implicit time discretization scheme which is
accepted custom in engineering and scientific computing, \cf\ \eg\
\cite{gajewski93}.

Please note that in device simulation one is always confronted with
contacted devices of heterogeneous material composition. That leads to
mixed boundary conditions and jumping material coefficients in the
model equations. Hence, standard theorems on existence, uniqueness
and regularity do not apply.

\section[Van Roosbroeck's system]{Van Roosbroeck's system}

\label{sec:setting}

\subsubsection*{Basic variables}
In the following we investigate van Roosbroeck's model for a
semiconductor device which describes the flow of electrons and holes
in a self-consistent electrical field due to drift and diffusion. The
physical quantities one is interested in are:
                                %
                                % \begin{itemize}
                                %                                 %
                                % \item 
the densities $u_1$ and $u_2$ of holes and electrons, 
                                %                                 %
                                % \item 
the densities $j_1$ and $j_2$ of the hole and electron current,
                                %                                 %
                                % \item 
the electrostatic potential $\Varphi$ of the self-consistent
electrical field, and
                                %                                 %
                                % \item 
the electrochemical potentials $\Fullphi_1$ and $\Fullphi_2$ of holes
and electrons
                                %                                 %
                                % \end{itemize}
                                %
These unknowns have to satisfy Poisson's equation and the
current--continuity equations for electrons and holes with some side
conditions. The latter are given by the relations between the
potentials and the densities.

\subsubsection*{Spatial domain}
We study only semiconductor devices which are quasi translational
invariant in one space direction or angular symmetric. In that case
van Roosbroeck's system in real space can be reduced to a similar set
of equations in the plane.  That means, we regard a cut through the
device perpendicular to the direction of invariance. Let
$\widehat{\Omega}$ be the resulting two-dimensional (bounded)
representative domain.  Parts of the device may be insulating, \eg\ 
formed by an oxide.  Then, electrons and holes can move only in a
sub-domain $\Omega$ of $\widehat{\Omega}$. This also covers the case
of charges which are artificially immobilized on a sub-domain
$\widehat{\Omega}\setminus\Omega$.
Furthermore, we mark out a part $\widehatNeumann$ of the boundary of
$\widehat{\Omega}$ where the device borders on an insulator. The
remaining part 
                                %$\widehatDirichlet$ 
of the boundary represents (possibly several) contacts of the device.
We also mark out a part $\Neumann$ of $\Omega$'s boundary. In the case
of a stand alone drift--diffusion model of the semiconductor device
again $\Neumann$ represents areas of the device bordering to an
insulator, whereas the remaining part
%% $\Dirichlet\df\interior{\partial\Omega\setminus\Neumann}$ 
is the contact area.

\subsubsection*{External control} 

In real--world modeling of semiconductor devices van Roosbroeck's
system often serves as a component in a compound model of the device.
Then the superordinated system --- for instance a circuit model ---
may exercise a control on van Roosbroeck's system.
Apart of a superordinated circuit model, compound models comprising in
addition to van Roosbroeck's system equations for the lattice
temperature or the power of lasing modes play an important role in
device simulation, \cf\ \eg\ 
\cite{gajewski93,spie-mqw:00,kare:smqw,bhk:02}.
But the concept of external control also comes to bear in segmentation
of the simulation domain, in particular in connection with multiscale
modeling, \cf\ \eg\ \cite{kare97b,kare97c,kare:current}.

If van Roosbroeck's equations serve as a component of a compound
model, then system parameters, state equations, boundary conditions,
et alii, possibly bear a different physical meaning than in the
stand-alone model.

We make assumptions about an external control from
the initial time $\TNULL$ up to a time $\TEINS$.

\subsection{Poisson equation}
The solution of the Poisson equation with mixed boundary conditions,
\begin{equation}
  \label{eq:poi}
  \begin{aligned}[2]
    -\dive 
    \left(
      \varepsilon \grad \Varphi
    \right) 
    %%  = \doping(t) + \extension u_1 - \extension u_2 
    & = \doping(t) + u_1 - u_2 
    &\qquad
    &\text{on $\widehat{\Omega}$,}
    \\
    \Varphi 
    &= \Dbvarphi(t) 
    &\qquad 
    &\text{on 
      \begin{math}
        \widehatDirichlet
        \df
        \interior{\partial\widehat{\Omega}\setminus\widehatNeumann}
      \end{math},}
    \\
    {\nu}\cdot{ 
      \left(
        \varepsilon \grad \Varphi
      \right)
    }
    + 
    \varepsilon_{\widehatNeumann} \Varphi 
    & =\Rbvarphi(t) 
    &\qquad
    &\text{on ${\widehatNeumann}$,}
  \end{aligned}
\end{equation}
gives the electrostatic potential $\Varphi$ on $\widehat{\Omega}$
subject to the electron and hole density $u_2$ and $u_1$.  Strictly
speaking, the densities $u_k$, $k=1,2$, are only defined on $\Omega$
but, we extend them by zero to $\widehat{\Omega}$.

The parameters in \eqref{eq:poi} have the following meaning:
$\varepsilon$ is a bounded, measurable function on $\widehat{\Omega}$
with values in the set of real, symmetric, $2\times2$, positive
definite matrices and corresponds to the spatially varying dielectric
permittivity on the space region occupied by the device.  Moreover, we
assume
\begin{equation*}
  \norm{\varepsilon(x)}_{\mathcal{B}(\mathbb{R}^2)}
  \le 
  \upp{\varepsilon} 
  \; \text{and} \;
  {(\varepsilon(x)\xi)}\cdot{\xi}
  \ge 
  \low{\varepsilon}
  \norm{\xi}_{\mathbb{R}^2}^2
  \quad
  \text{for almost all $x\in\widehat{\Omega}$ 
    and all $\xi\in\mathbb{R}^2$}
\end{equation*}
with two strictly positive constants $\low{\varepsilon}$ and
$\upp{\varepsilon}$.  Furthermore, $\varepsilon_{\widehatNeumann}$ is
a non-negative function on ${\widehatNeumann}$, representing the
capacity of the part of the device surface bordering on an insulator.
We assume that $\widehatDirichlet$ is not empty or
$\varepsilon_{\widehatNeumann}$ is positive on a subset of
$\widehatNeumann$ with positive arc measure.  In other words, the
device has a Dirichlet contact or part of its surface has a positive
capacity.
$\Dbvarphi(t)$ and $\Rbvarphi(t)$ are the voltages applied at the
contacts of the device, and $\doping(t)$ represents a charge.
In the case of a stand alone drift--diffusion model $\Dbvarphi$,
$\Rbvarphi$, and $\doping$ are constant in time, and
$\doping$ solely is the charge density of dopants in the
semiconductor materials composing the device.
In general, $\Dbvarphi$, $\Rbvarphi$, and $\doping$ are function
which are defined on the time interval $[\TNULL,\TEINS]$ where a possible
control acts on the device.

\subsection{Current--continuity equations}                         

The current--continuity equations for holes and electrons ($k=1,2$,
respectively)
\begin{equation}
  \label{eq:cc}
  u_k' 
  - \dive j_k 
  =
  \reaction_k(t,\Varphi,\Fullphi_1,\Fullphi_2)
  \qquad \text{on $\Omega$}
\end{equation}
characterize the evolution of the electron and hole density under the
action of the currents $j_k$ and the reactions $\reaction_k$ subject to
the mixed boundary conditions
\begin{equation}
  \label{eq:cc:bc}
  \begin{aligned}[2]
    \Fullphi_k(t) 
    &=\Bvphi{k}(t) 
    &\qquad 
    &\text{on $\Dirichlet\df\interior{\partial\Omega\setminus\Neumann}$,} 
    \\
    {\nu}\cdot{j_k} 
    &= 0 
    &\qquad 
    &\text{on $\Neumann$,} 
  \end{aligned}
                                %  \qquad k=1,2,
                                %  \quad t\in[\TNULL,\TEINS],
\end{equation}
from the initial conditions
\begin{equation}
  \label{eq:cc:ini}
  \Fullphi_k(\TNULL)=\iniphi_k. 
                                %  \qquad k=1,2.
\end{equation}

Each $\reaction_k$, $k=1,2$ is a reaction term which models the
generation and annihilation of electrons and holes.  In particular,
this term covers the recombination of electrons and holes in the
semiconductor device. $\reaction_1$ and $\reaction_2$ can be rather
general functions of the particle and current densities, \cf\ 
\secref{sec:reactions}.  We require that the set
$\Dirichlet=\interior{\partial\Omega\setminus\Neumann}$ is not empty.
The boundary values $\Bvphi{1}$, $\Bvphi{2}$ in general depend on
time. Moreover, the reactions $\reaction_k$ may explicitly depend on
time. This dependence on time, again, allows for a control of the
system by some other part of a superordinated compound model.

\subsection{Carrier and current densities}

Van Roosbroeck's system has to be complemented by a prescription
relating the density of electrons and holes as well as the densities
of the electron and hole current to the chemical potentials of these
charge carriers. We assume
\begin{equation}
  \label{eq:densities}
  u_k(t,x) 
  \df
  \rho_k(t,x) 
  \mathcal{F}_k 
  \left(
    \chempot_k(t,x)
  \right)\, ,
   \quad  x \in \Omega 
  ,
  \qquad k=1,2,
\end{equation}
where $\chempot_1$ and $\chempot_2$ are the chemical potentials
\begin{equation}
  \label{eq:chempot}
  \chempot_k 
  \df
  \Fullphi_k
  + (-1)^k \Varphi 
  + b_k
%  , \qquad x\in{\Omega}
  , \qquad k=1,2,
\end{equation}
and $\Fullphi_2$, $\Fullphi_1$ are the electrochemical potentials of
electrons and holes, respectively.
$b_k$, $\rho_k$, $k=1,2$ are positive, bounded functions on $\Omega$.
They describe the electronic properties of the materials composing the
device. $b_2$ and $b_1$ are the band edge offsets for electrons and
holes, and $\rho_2$, $\rho_1$ are the corresponding effective band
edge densities of states.
If the equations under consideration form part of a compound model for
the semiconductor device, then $b_k$, $\rho_k$, $k=1,2$, may depend on
time. For instance, the $\rho_k$ could be subject to an external
control of the device temperature. Then they depend on time via the
temperature. Mathematically, we assume the  following. 
\begin{assu}
  \label{assu:effband}
  For every $t\in[\TNULL,\TEINS]$ the functions $\rho_k(t)$ are
  essentially bounded on $\Omega$ and admit positive lower bounds
  which are uniform in $t\in[\TNULL,\TEINS]$.  The mappings
  \begin{equation} 
    \label{e-rho1}
    [\TNULL,\TEINS] \ni t \mapsto \rho_k(t) \in L^2(\Omega), \quad k=1,2
  \end{equation}
  are differentiable on the interval $]\TNULL,\TEINS[$ with H\"older
  continuous derivatives $\rho_k'$.
\end{assu}
The functions $\mathcal{F}_1$ and $\mathcal{F}_2$ represent the
statistical distribution of the holes and electrons on the energy
band.  In general, Fermi--Dirac statistics applies, \ie\
\begin{equation}
  \label{eq:fermi-integral}
  \mathcal{F}_k (s) 
  \df 
  \frac{2}{\sqrt{\pi}}  
  \int_0^\infty 
  \frac{\sqrt{t}}{1+\exponential{t-s}}
  \der{t}
  ,
  \qquad
  s\in\mathbb{R}
  .
\end{equation}
However, often Boltzmann statistics 
\begin{math}
  \label{eq:boltzmann}
  \mathcal{F}_k (s) 
  = \exponential{s}
\end{math}
is a good approximation.
                                
As for the kinetic relations specifying the current--continuity
equations we assume that the electron and hole current is driven by
the negative gradient of the electrochemical potential of electrons
and holes, respectively. More precisely, the
current densities are given by
\begin{equation}
  \label{eq:curr-dens}
  j_k(t,x) 
  = 
  -
  \mathcal{G}_k 
  \left(
    \chempot_k (t,x)
  \right)
  \mu_k(x) 
  \,
  \grad \Fullphi_k(t,x)\; , 
  \quad x \in \Omega
  ,
  \qquad k = 1,2. 
\end{equation}
The mobilities $\mu_2$ and $\mu_1$ for the electrons and holes,
respectively, are measurable, bounded function on ${\Omega}$ with
values in the set of real, $2\times2$, positive definite matrices
satisfying for almost all $x\in\widehat{\Omega}$ and all
$\xi\in\mathbb{R}^2$
\begin{equation*}
  \norm{\mu_k(x)}_{\mathcal{B}(\mathbb{R}^2)}
  \le \upp{\mu}
  \quad \text{and} \quad
  {(\mu_k(x)\xi)}\cdot{\xi}
  \ge 
  \low{\mu}
  \norm{\xi}_{\mathbb{R}^2}^2
  ,
  \qquad 
  k=1,2,
\end{equation*}
with two strictly positive constants $\low{\mu}$ and $\upp{\mu}$.
The mobilities are accounted for on the parts of the device where
electrons and holes can move due to drift and diffusion.
\begin{rem}
  \label{rem:k}
  In semiconductor device modeling, usually, the functions
  $\mathcal{G}_k$ and $\mathcal{F}_k$ coincide, \cf\ \eg\
  \cite{selberherr84} and the references there. However, a rigorous
  formulation as a minimal problem for the free energy reveals that
  \begin{math}
    \mathcal{G}_k = \mathcal{F}_k^\prime
  \end{math}
  is appropriate. This topic has been thoroughly investigated for
  analogous phase separation problems, \cf\ 
  \cite{quastel:92,quastel:99,lebowitz:97,lebowitz:98}, \cf\ also
  \cite{skrypnik:02} and \cite{griepentrog04}.  In order to cover both
  cases we regard independent functions $\mathcal{G}_k$ and
  $\mathcal{F}_k$.
\end{rem}
\begin{assu}
  \label{assu:distri}
  Mathematically, we demand that the distribution functions
  $\mathcal{F}_k$, $\mathcal{G}_k$, $k=1,2$, are defined on the real
  line, take positive values, and are
                                %$\from \mathbb{R} \to (0,\infty)$
  either exponentials, or twice continuously differentiable and
  polynomially bounded. Moreover, $\mathcal{F}_1'$, $\mathcal{F}_2'$
  are strictly positive on $\mathbb{R}$. In the sequel we will call
  such distribution functions 'admissible.'  This includes Boltzmann
  statistics, as well as Fermi--Dirac statistics (\cf\ 
  \eqref{eq:fermi-integral}).
\end{assu}
Let us comment on the (effective) band edges $b_k$ and the (effective)
densities of states $\rho_k$, \cf\ \eqref{eq:densities} and
\eqref{eq:chempot}:
%\begin{rem}
%  \label{rem:bands}
Basically the band edge offsets $b_k$ and the effective band edge
densities of states $\rho_k$ are material parameters. In a
heterogeneous semiconductor device they are generically piecewise
constant on the spatial domain $\Omega$.
As Assumption~\ref{assu:bands} reveals, we cannot cope with such a
situation as far as the band edges $b_k$ are concerned. However, in
the case of Boltzmann statistics one can rewrite \eqref{eq:densities}
and \eqref{eq:chempot} as
\begin{equation*}
  u_k 
  =
  \rho_k
  \exponential{b_k}
  \exponential{\left(
      \Fullphi_k
      + (-1)^k \Varphi 
    \right)}
  \; \;
  \text{on} \; {\Omega}
  ,
  \qquad k=1,2,
\end{equation*}
with modified effective densities of states and
identically vanishing band edge offsets.
In the case of Fermi--Dirac statistics this reformulation is not
possible and one has to recourse to some approximation of the $b_k$
by functions confirming to Assumption~\ref{assu:bands}.
Discontinuities of the band edge offsets up to now seem to be an
obstacle in whatever approach to solutions of van
Roosbroeck's equations, if the statistical distribution function is
not an exponential, \cf\ \eg\ \cite{pub:938}.

There are compound multiscale models of semiconductor devices such
that the effective band edges and the effective densities of states
result by upscaling from quantum mechanical models for the electronic
structure in heterogeneous semiconductor materials, \cf\ 
\cite{spie-mqw:00,bhk:02,pub:1133}. In view of an offline coupling to
electronic structure calculations we allow for an explicit dependence
of $\rho_k$, and $b_k$ on time.

%% We allow for an explicit
%% dependence of $\rho_k$, and $b_k$ on time in view of this and of offline
%% coupling to electronic structure calculations.

\subsection{Reaction rates}
\label{sec:reactions}

The reaction terms on the right hand side of the current--continuity
equations can be rather general functions of time, of the
electrostatic potential, and of the vector of the electrochemical
potentials. $\reaction_1$ and $\reaction_2$ describes the production
of holes and electrons, respectively --- generation or annihilation,
depending on the sign of the reaction term. Usually van Roosbroeck's
system comprises only recombination of electrons and holes:
\begin{math}
  \reaction=\reaction_1=\reaction_2.
\end{math}
We have formulated the equations in a  more general way, in
order to include also coupling terms to other equations of a
superordinated compound model. That is why we also allow for an
explicit time dependency of the reaction rates.

Our formulation of the reaction rates, in particular, includes a
variety of models for the recombination and generation of
electrons--hole pairs in semiconductors.  This covers non-radiative
recombination of electrons and holes like the Shockley--Read--Hall
recombination due to phonon transition and Auger recombination.
% (three particle transition). 
But, radiative recombination (photon
transition), both spontaneous and stimulated, is also included.
Mathematical models for stimulated optical recombination typically
require the solution of additional equations for the optical
field. Thus, the recombination rate may be a non-local operator.
Moreover, by coupling van--Roosbroecks system to the optical field
some additional control of this optical field may also interact with
the internal electronics. For instance, in modeling and simulation of
edge--emitting multiple--quantum--well lasers van--Roosbroeck's system
augmented by some Helmholtz equation often serves as a transversal (to
the light beam) model, and a control of the optical field is exercised
by a master equation or some model for the longitudinal (on the axis
of the light beam) behavior of the laser, \cf\ \eg\
\cite{wuensche:etal93,spie-mqw:00,bhk:02}.

Modeling recombination of electron--hole pairs in semiconductor
material is an art in itself, \cf\ \eg\ \cite{landsberg91}.
However, for illustration, let us list some common recombination
models, \cf\ \eg\ \cite{selberherr84,gajewski93} and the references
cited there.

\emph{Shockley--Read--Hall recombination} 
(phonon transitions):
\begin{equation*}
  \reaction_1
  =\reaction_2
  =\reaction^{\mathrm{SRH}} 
  =
  \frac{u_1 u_2 - n_i^2}{\tau_2(u_1+n_1)+\tau_1(u_2+n_2)},
\end{equation*}
where $n_i$ is the intrinsic carrier density, $n_1$, $n_2$ are
reference densities, and $\tau_1$, $\tau_2$ are the lifetimes of
holes and electrons, respectively.  $n_i$, $n_1$, $n_2$, and $\tau_1$,
$\tau_2$ are parameters of the semiconductor material; thus, depend
on the space variable, and ultimately, also on time.

\emph{Auger recombination} (three particle transitions):
\begin{equation*}
  \reaction_1
  =\reaction_2
  =\reaction^{\mathrm{Auger}} 
  =
  (u_1 u_2 - n_i^2)
  (c_1^{\mathrm{Auger}} u_1 + c_2^{\mathrm{Auger}}u_2),
\end{equation*}
where $c_1^{\mathrm{Auger}}$ and $c_2^{\mathrm{Auger}}$ are the Auger
capture coefficients of holes and electrons, respectively, in the
semiconductor material.

\emph{Stimulated optical recombination:} 
\begin{equation*}
  \reaction_1
  =\reaction_2
  =\reaction^{\mathrm{stim}}
  =
  \sum_j 
  f(\sigma_j)
  \frac{\abs{\psi_j}^2}{\int\abs{\psi_j}^2},
\end{equation*}
where $f$ additionally depends on the vector of the densities, and on
the vector of the electrochemical potentials.  $\sigma_j$, $\psi_j$
are the eigenpairs of a scalar Helmholtz--operator:
\begin{equation*}
  \Delta \psi_j + \epsilon(u_1,u_2) \psi_j = \sigma_j \psi_j
  .
\end{equation*}
                                %
                                %with Sommerfeld's radiation condition. 
In laser modeling each eigenpair corresponds to an optical (TE) mode of the
laser and $\abs{\psi_j}^2$ is the intensity of the electrical field of
the $\sigma_j$--mode. $\epsilon$ is the dielectric permittivity (for
the optical field); it depends on the density of electrons and holes.
The scalar Helmholtz--equation originates from the Maxwell equations
for the optical field \cite{wuensche91}. 

The functional analytic requirements on the reaction terms will be
established in Assumption~\ref{assu:recomb}.

\section[Mathematical prerequisites]{Mathematical prerequisites}
\label{sec:pre}

In this section we introduce some mathematical terminology and make
precise assumptions about the problem.

\subsection{General Assumptions}
\label{sec:general}

For a Banach space $X$ we denote its norm by $\norm{\cdot}_X$ and the
value of a bounded linear functional $\psi^*$ on $X$ in $\psi\in{X}$
by $\dual{\psi^*}{\psi}_{X}$. If $X$ is a Hilbert space, identified
with its dual, then $\dual{\cdot}{\cdot}_{X}$ is the scalar product in
$X$.  Just in case $X$ is the space $\mathbb{R}^2$, the scalar
product of $a,b\in\mathbb{R}^2$ is written as ${a}\cdot{b}$.
Upright $\xoplus{X}$ denotes the direct sum $X{\oplus}X$ of slanted
$X$ with itself.
$\mathcal{B}(X;Y)$ is the space of linear, bounded operators from $X$
into $Y$, where $X$ and $Y$ are Banach spaces. We abbreviate
$\mathcal{B}(X)=\mathcal{B}(X;X)$ and we denote by
$\mathcal{B}_\infty(X)$ the space of linear, compact operators on the
Banach space $X$.
The notation $[X,Y]_\theta$ means the complex interpolation space of
$X$ and $Y$ to the index $\theta \in[0,1]$.
The (distributional) $\nabla$--calculus applies. If $\psi$ is a
(differentiable) function on an interval taking its values in a Banach
space, then $\psi'$ always indicates its derivative.

\subsection{Spatial Domains}
\label{sec:domain}

Throughout this paper we assume that $\widehat{\Omega}$ as well as
$\Omega$ are bounded Lipschitz domains in $\mathbb{R}^2$, \cf\ 
\cite[Ch.~1]{grisvard85}.  By $\extension$ we denote the operator
which extends any function defined on $\Omega$ by zero to a function
defined on $\widehat{\Omega}$.  Conversely, $\restrict$ denotes the
operator which restricts any function defined on $\widehat{\Omega}$ to
$\Omega$. The operators $\extension$ and $\restrict$ are adjoint to
each other with respect to the duality induced by the usual scalar
product in spaces of square integrable functions.

With respect to the marked out Neumann boundary parts
$\widehatNeumann\subset\partial\widehat{\Omega}$ and
$\Neumann\subset\partial\Omega$ of the boundary of $\widehat{\Omega}$
and $\Omega$ we assume each being the union of a finite
set of open arc pieces such that no connected component
                                %Zusammenhangskomponente
of $\partial\widehat{\Omega}\setminus{\widehatNeumann}$ and
$\partial\Omega\setminus\Neumann$ consists only of a single point. We
denote the parts of the boundary where Dirichlet boundary conditions
are imposed by
\begin{math}
  \widehatDirichlet
  \df
  \interior{\partial\widehat{\Omega}\setminus\widehatNeumann}
\end{math}
and
\begin{math}
  \Dirichlet\df\interior{\partial\Omega\setminus\Neumann}.
\end{math}

\subsection{Function spaces and linear elliptic operators}
\label{sec:spaces}

We exemplarily define spaces of real-valued functions on spatial
domains with respect to the bounded domain $\Omega\subset\mathbb{R}^2$
and its boundary. Spaces of functions on $\widehat{\Omega}$ and parts
of its boundary may be similarly defined and are denoted by hatted
symbols.

If $r\in[1,\infty[$, then $L^r$ is the space of real, Lebesgue
measurable, $r$-integrable functions on $\Omega$ and $L^\infty$ is the
space of real, Lebesgue measurable, essentially bounded functions on
$\Omega$.
$W^{1,r}$ is the usual Sobolev space
$W^{1,r}(\Omega)$,
%space of Bessel potentials with the
%differentiability index $s\in[0,1]$ and integrability index
%$r\in[1,\infty[$ on $\Omega$,
\cf\ \eg\ \cite{triebel}.
$W^{1,r}_{\Neumann}$ is the closure in $W^{1,r}$ of
\begin{equation*}
  \left\{
    \psi|_{\Omega}
    \with
    \psi \in C^\infty_0(\mathbb{R}^2)
    , \;
    \supp \psi 
    \cap
    (\partial\Omega \setminus \Neumann)
    =\emptyset 
  \right \}
  ,
\end{equation*}
\ie \, $W^{1,r}_{\Neumann}$ consists of all functions from
$W^{1,r}$ with vanishing trace on $\Dirichlet$.
$W^{-1,r}_{\Neumann}$ denotes the dual of
$W^{1,r'}_{\Neumann}$, where $\textfrac{1}{r}+\textfrac{1}{r'}=1$.
\begin{math}
  \dual{\cdot}{\cdot}_{W^{1,2}_{\Neumann}}
\end{math}
is the dual pairing between $W^{1,2}_{\Neumann}$ and $W^{-1,2}_{\Neumann}$.
Correspondingly, the divergence for a vector of square integrable
functions is defined in the following way: If $j\in\xoplus{L}^2$, then
\begin{math}
  \dive j \in W^{-1,2}_{\Neumann}
\end{math}
is given by
\begin{equation}
  \label{eq:div}
  \lrdual{\dive j}{\psi}_{W^{1,2}_{\Neumann}}  
  = 
  -\int_\Omega 
  {j}\cdot{\grad \psi} \der{x},
  \qquad 
  \psi \in W^{1,2}_{\Neumann}.
\end{equation}

$\sigma$ is the natural arc measure on the boundary of $\Omega$. 
We denote by
$L^\infty(\partial\Omega)$ and 
$L^r(\partial\Omega)$,
the spaces of $\sigma$-measurable, essentially bounded, and
$r$-integrable, $r\in[1,\infty[$, functions on $\partial\Omega$,
respectively. Moreover,
$W^{s,r}(\partial\Omega)$
denotes the Sobolev space of fractional order $s\in]0,1]$ and
integrability exponent $r\in[1,\infty[$ on $\partial\Omega$, \cf\ 
\cite[Ch.~1]{grisvard85}.
Mutatis mutandis for functions on ${\sigma}$-measurable, relatively
open parts of $\partial{\Omega}$.

Let us now define in a strict sense the (linear) Poisson operator and
the elliptic operators governing the current continuity equations.
\begin{defn} 
  \label{d-opera}
  We define the Poisson operator 
  \begin{math}
    -\dive\varepsilon\grad 
    \from{
      \widehat{W}^{1,2} 
      \to
      \widehat{W}_{\widehatNeumann}^{-1,2}
    }
  \end{math}
  by
  \begin{equation} 
    \label{e-poi}
    \dual{
      -\dive\varepsilon\grad\psi_1}
    {\psi_2}_{\widehat{W}_{\widehatNeumann}^{1,2}}
    \df  
    \int_{\widehat{\Omega}}
    {\varepsilon \grad\psi_1}\cdot{\grad \psi_2}
    \der{x}
    +
    \int_{\widehatNeumann}
    \varepsilon_{\widehatNeumann} \psi_1 \psi_2 
    \der{\widehat{\sigma}}
    ,
%%     \\
%%     \psi_1 \in \widehat{W}^{1,2},
%%     \quad
%%     \psi_2 \in \widehat{W}_{\widehatNeumann}^{1,2}.
  \end{equation}
  for 
  $\psi_1\in\widehat{W}^{1,2}$ and
  $\psi_2\in\widehat{W}_{\widehatNeumann}^{1,2}$.
  $\Linpoisson$ denotes
  the restriction of $-\nabla\cdot\varepsilon\nabla$ to
  $\widehat{W}_{\widehatNeumann}^{1,2}$; we denote the maximal
  restriction of $\Linpoisson$ to any range space which continuously
  embeds into $\widehat{W}_{\widehatNeumann}^{-1,2}$ by the same
  symbol $\Linpoisson$.
\end{defn}
\begin{defn}
  \label{def:cont}
  With respect to a function $\varsigma\in{}L^\infty$
  we define the operators
  \begin{multline*}
    -\dive{\varsigma\mu_k}\grad
    \from{W^{1,2} \to W^{-1,2}_{\Neumann}}
    ,
    \quad
    k=1,2,
    \quad
    \text{by}
    \\
    \dual{-\dive \varsigma\mu_k \grad\psi_1}{\psi_2}_{W^{1,2}_{\Neumann}}
    \df
    \int_\Omega
    \varsigma\;
    {\mu_k\grad\psi_1}\cdot{\grad\psi_2}
    \der{x}
    ,
    \qquad
    \psi_1\in W^{1,2},
    \;
    \psi_2\in W^{1,2}_{\Neumann}
    .
  \end{multline*}
  If, in particular, $\varsigma\equiv{1}$, then we simply write
  $\fulla{k}$ for $-\dive\mu_k\grad$.  Moreover, we denote the
  restriction of $\fulla{k}$ to the space $W^{1,2}_{\Neumann}$ by
  $a_k$, \ie\ 
  \begin{math}
    a_k \from{W^{1,2}_{\Neumann} \to W^{-1,2}_{\Neumann}}.
  \end{math}
                                %
%  Finally, we define
%                                %
%  \begin{equation*}
%    \fulla
%    \df
%    \left(
%      \begin{smallmatrix}
%        \fulla{1} & 0
%        \\
%        0         & \fulla{2}
%      \end{smallmatrix}
%    \right)
%    \from{\xoplus{W}^{1,2} \to \xoplus{W}^{-1,2}_{\Neumann}}
%    \quad \text{and}\quad
%    a
%    \df
%    \left(
%      \begin{smallmatrix}
%        a_1 & 0
%        \\
%        0   & a_2
%      \end{smallmatrix}
%    \right)
%    \from{\xoplus{W}^{1,2}_{\Neumann} \to \xoplus{W}^{-1,2}_{\Neumann}}
%    .
%  \end{equation*}
%                                %
\end{defn}
\begin{prop}
  \label{prop:isomorphy}
  \emph{(\cf\ \cite{groeger89} and \cite{groeger:rehberg89})} There is
  a number $\hat{q}>2$ (depending on $\widehat{\Omega}$, $\varepsilon$
  and $\widehatNeumann$) such that for all $q\in[2,\hat{q}]$ the
  operator
  \begin{math}
    \Linpoisson 
    \from{
      \widehat{W}^{1,q}_{\widehatNeumann} 
      \to
      \widehat{W}^{-1,q}_{\widehatNeumann}
    }
  \end{math}
  is a topological isomorphism.  Moreover, there is a $\check{q}>2$
  (depending on $\Omega$, $\mu_1$, $\mu_2$ and $\Neumann$) such that
  for all $q\in[2,\check{q}]$ the operators
  \begin{math}
    a_k 
    \from{
      W^{1,q}_{\Neumann} 
      \to
      W^{-1,q}_{\Neumann}
    }  
  \end{math}
  provide topological isomorphisms, and additionally, generate analytic
  semigroups on $W^{-1,q}_{\Neumann}$.
\end{prop}
\begin{defn} 
  \label{def:pq}
  From now on we fix a number $q\in]2,\min(4,\hat{q},\check{q})[$ and
  define $p\df\frac{q}{2}$.
                                %  
%%   Furthermore, we fix for all what follows a number 
%%   $\theta\in]\frac{1}{2}+\frac{2}{q},1]$.
                                %
  With respect to this $p$ we define the operators
  \begin{gather*}
    A_k \from{\psi \mapsto a_k \psi}
    ,\quad
    \psi \in \domainA_{k} 
    \df \dom(A_k) 
    \df
    \left\{
      \psi \in W^{1,2}_{\Neumann}
      \with
      a_k\psi \in L^p
    \right\}
    ,\quad
    k=1,2,
    \\
    \matrixA \from{ \domainA \to \xoplus{L}^p }
    ,\quad
    \matrixA \df 
    \left(
      \begin{smallmatrix}
        A_1 & 0
        \\
        0   & A_2
      \end{smallmatrix}
    \right)
    ,\quad
    \domainA \df \dom(\matrixA) = \domainA_{1} \oplus \domainA_{2}
    \embedcontinuously 
    \xoplus{L}^p
    .
  \end{gather*}
\end{defn}
\begin{rem}
  \label{rem:normnull}
  If $\psi\in\domainA_{k}$, $k=1,2$, then  
  \begin{math}
    {\nu}\cdot{(\mu_k \grad \psi)}
    |_{\Neumann} 
    =
    0
  \end{math}
  in the sense of distributions, \cf\ \eg\
  \cite[Ch.~1.2]{ciarlet79} or
  \cite[Ch.1.2]{gajewski:groeger:etal74}.
\end{rem}
After having fixed the number $q$ and, correspondingly, the space
$L^p$, we will now formulate our mathematical requirements on the
reaction terms:
\begin{assu}
  \label{assu:recomb}
  The reaction terms $\reaction_k$, $k=1,2$,  are mappings
  \begin{equation*}
    \reaction_k \from{[\TNULL,\TEINS] 
      \times \widehat{W}^{1,q}
      \times \xoplus{W}^{1,q}
      \to
      L^p
    }.
  \end{equation*}
  Moreover, we assume that there is a real number $\eta\in]0,1]$ and
  for any bounded subset
  $M\subset\widehat{W}^{1,q}\oplus\xoplus{W}^{1,q}$ a constant
  $\reaction_M$ such that
  \begin{multline*}
    \lrnorm{
      \reaction_k(\ta,v,\psi)
      - 
      \reaction_k(\tb,\check{v},\check{\psi})
    }_{L^p} 
    \\
    \le
    \reaction_M
    \left(
      \abs{\ta-\tb}^\eta
      +
      \norm{
        v - \check{v}
      }_{\widehat{W}^{1,q}} 
      +
      \norm{
        \psi - \check{\psi}
      }_{\xoplus{W}^{1,q}}
    \right)
    ,
    \\
    \ta,\tb \in [\TNULL,\TEINS]
    ,
    \quad
    (v,\psi),\, (\check{v},\check{\psi}) \in M
    .
  \end{multline*}
\end{assu}
\begin{assu}
  \label{assu:bands}
  The functions
  \begin{math}
    b_k
    \from{
      [\TNULL,\TEINS]
      \to 
      W^{1,q}
    },
  \end{math}
  $k=1,2$, are H\"older continuous.  Moreover, they are H\"older
  continuously differentiable when considered as $L^p$ valued.
\end{assu}

\subsection{Representation of Dirichlet boundary values}
\label{sec:diribv}

For setting up the Poisson and current--continuity equations in
appropriate function spaces we must split up the solution into parts,
where one part represents the inhomogeneous Dirichlet boundary values
$\Dbvarphi$ and $\Bvphi{k}$, $k=1,2$. In this section we
treat of just this representation.
We make the following assumptions about the Dirichlet boundary values
of the electrochemical potentials $\phi_k$, $k=1,2$, and for their
initial values, \cf\ \eqref{eq:cc:bc}, \eqref{eq:cc:ini}.
\begin{assu}
  \label{assu:cc:bc} 
  There is a H\"older continuous function
  \begin{equation*}
    \bvphi{}=(\bvphi{1},\bvphi{2}) \from{ [\TNULL,\TEINS] \to \xoplus{W}^{1,q}},
    \quad k=1,2,
  \end{equation*}
                                %$k=1,2$, 
  such that for all $t\in[\TNULL,\TEINS]$
  \begin{align}
    \label{eq:cc:bv1}
    \fulla{k} \bvphi{k}(t) 
    & = 0
    \\ 
    \label{eq:cc:bv2}
    \trace 
    \big( \bvphi{k}(t) \big)
    \big|_{\Dirichlet} &=\Bvphi{k}(t) 
  \end{align}
  Moreover, we assume, that each $\bvphi{k}$, $k=1,2$,
  --- as a function with values in $L^p$ --- 
  is differentiable and its derivative is H\"older continuous.
\end{assu}
\begin{rem} 
  \label{rem:curnormal}
  It should be noted that \eqref{eq:cc:bv1} and the definition of the
  operators $\fulla{k}$ imply $\nu\cdot\mu_k\grad\bvphi{k}=0$ on
  $\Neumann$ in the distributional sense, \cf\ \eg\ 
  \cite[Ch.~1.2]{ciarlet79} or
  \cite[Ch.~II.2]{gajewski:groeger:etal74}.  This implies for the
  current densities \eqref{eq:curr-dens} that $\nu\cdot{}j_k=0$ on
  $\Neumann$ in the distributional sense, provided that 
  $\chempot_k\in{}W^{1,q}$.
\end{rem}
We will now give a sufficient condition on $\Bvphi{k}$ for the
existence of a $\bvphi{k}$ with the assumed properties.
\begin{lem} 
  \label{lem:diri}
  1.~If
  \begin{math}
    \psi \in W^{1-\textfrac{1}{q},q}(\Dirichlet),
  \end{math}
  then there is a unique function $\bvpsi\in{W^{1,q}}$ fulfilling 
  \begin{equation*}
    \fulla{k} \bvpsi
    = 0,
    \quad\text{and}\quad
    \trace 
    (\bvpsi)
    \big|_{\Dirichlet} 
    = \psi. 
  \end{equation*}
  \par
  2.~If
  \begin{math}
    \psi 
    \from{ 
      [\TNULL,\TEINS] \to
      W^{1-\textfrac{1}{q},q}(\Dirichlet)
    }
  \end{math}
  is H\"older continuous with index $\eta$, then the function
  \begin{math}
    \bvpsi 
    \from{ [\TNULL,\TEINS] \to W^{1,q} }
  \end{math}
  which is given for each $t\in[\TNULL,\TEINS]$ by item~1 is also H\"older
  continuous with index $\eta$.
  Moreover, if $\psi$ --- as a function with values in
  $W^{\textfrac{1}{2},2}(\Dirichlet)$ --- is H\"older continuously
  differentiable with H\"older index $\eta$, then $\bvpsi$ is H\"older
  continuously differentiable with H\"older index $\eta$.
\end{lem}
\begin{proof}
  Let 
  \begin{math}
    \operatorname{ex}
    \from{
      W^{1-\textfrac{1}{q},q}(\Dirichlet)
      \to
      W^{1-\textfrac{1}{q},q}(\partial\Omega)
    }
  \end{math}
  be a linear and continuous extension operator, and let $\trace^{-1}$
  be a linear and continuous right inverse of the trace operator
  \begin{math}
    \trace 
    \from{
      W^{1,q} (\Omega)
      \to 
      W^{1-\textfrac{1}{q},q}(\partial\Omega)
    }.
  \end{math}
  Such operators exist according to \cite[Thm~1.4.3.1]{grisvard85} and
  \cite[Thm~1.5.1.3]{grisvard85}, respectively.
  Thus,
  \begin{math}
    \trace^{-1} 
    \circ 
    \operatorname{ex}
    \psi \in W^{1,q}.
  \end{math}
  Moreover, let $\wildbvpsi$ be the solution of the differential
  equation
  \begin{equation}
    \label{eq:deiff}
    a_k \wildbvpsi
    = 
    \fulla{k} 
    \circ
    \trace^{-1}
    \circ
    \operatorname{ex}
    \psi
  \end{equation}
  in $W^{1,q}_{\Neumann}$.
  This solution exists and is unique because the right hand side of
  \eqref{eq:deiff} is from $W^{-1,q}_{\Neumann}$ and the operators
  $a_k$ are isomorphisms from $W^{1,q}_{\Neumann}$ onto
  $W^{-1,q}_{\Neumann}$.
  We now define 
  \begin{equation}
    \label{eq:DefinitionExtendedDirichletData}
    \bvpsi \df       
    \trace^{-1}
    \circ
    \operatorname{ex}
    \psi
    -
    \wildbvpsi
    .
  \end{equation}
  The asserted properties of $\bvpsi$ follow directly from the
  construction.
  
  The second assertion is proved by observing that all steps in the
  first part of the proof depend linearly on the datum.
\end{proof}
\begin{assu} 
  \label{assu:cc:ini} 
  We assume that the initial values $\iniphi_k$ belong to $W^{1,q}$,
  $k=1,2$. Moreover, there is a
  $\theta\in]\textfrac{1}{2}+\textfrac{1}{q},1[$ such that for each of
  the initial values $\iniphi_k$ the difference
  $\iniphi_k-\bvphi{k}(\TNULL)$ belongs to the complex interpolation space
  $[L^p,\domainA_{k}]_{\theta}$.
\end{assu}
\begin{rem} 
  \label{rem:initial} 
  For all $\theta\in]\textfrac{1}{2}+\textfrac{1}{q},1[$ the space
  $[L^p,\domainA_{k}]_{\theta}$ compactly embeds into
  $W^{1,q}_{\Neumann}\embedcontinuously{}L^\infty$, \cf\ 
  \cite[Thm.~5.2]{pub:765}.
\end{rem}
With respect to the inhomogeneous terms $\Dbvarphi$ and $\Rbvarphi$ in
the boundary conditions of Poisson's equation \eqref{eq:poi} we make
the following assumptions.
\begin{assu}
  \label{assu:poi-diri-bv}
  There is a H\"older continuous function
  \begin{math}
    \dbvarphi\from{
      [\TNULL,\TEINS] 
      \to 
      \widehat{W}^{1,q}
    }
  \end{math}
  such that $\dbvarphi$ --- as a function from $[\TNULL,\TEINS]$ into
  $\widehat{L}^{p}$ --- is H\"older continuously
  differentiable. For all $t\in[\TNULL,\TEINS]$ it holds true
  \begin{align}
    \label{eq:trace1}
    - \nabla \cdot \varepsilon \nabla \dbvarphi(t) 
    & = 0, 
    \\
    \label{eq:trace2}
    \trace
    \big( \dbvarphi(t) \big)
    \big|_{\widehatDirichlet} 
    &=  \Dbvarphi(t).
  \end{align}
  The function 
  \begin{equation*}
    [\TNULL,\TEINS] \ni t 
    \mapsto 
    \Rbvarphi(t) \in L^\infty(\widehatNeumann)    
  \end{equation*}
  is differentiable and possesses a H\"older continuous derivative.
\end{assu}
\begin{rem} 
  \label{r-represent}
  Similar to \lemref{lem:diri} it is possible to give a sufficient
  condition on the existence of a representing function
  $t\mapsto\dbvarphi(t)$ which only rests on the function
  $t\mapsto\Dbvarphi(t)$. We do not carry out this here.
\end{rem}
\begin{rem} 
  \label{r-extend}
  For all $t\in[\TNULL,\TEINS]$ we extend $\Rbvarphi(t)$ by zero to a
  $\widehat{\sigma}$--measurable, essentially bounded function on
  $\partial\widehat{\Omega}$. Due to the continuous embedding
  \begin{displaymath}
    \widehat{W}_{\widehatNeumann}^{1,q'}
    \embedcontinuously 
    \widehat{W}^{1,q'}
    \embedcontinuously 
    W^{1-\textfrac{1}{q'},q'}(\partial\widehat{\Omega}) 
    \embedcontinuously 
    L^{q'}(\partial\widehat{\Omega})
    ,
  \end{displaymath}
  \cf\ \cite[Thm~1.5.1.3]{grisvard85}, there is a continuous embedding
  \begin{displaymath}
    L^{\infty}(\partial\widehat{\Omega})
    \embedcontinuously 
    L^{q}(\partial\widehat{\Omega})    
    \embedcontinuously 
    \widehat{W}_{\widehatNeumann}^{-1,q}
    .
  \end{displaymath}
  Thus, $\Rbvarphi(t)$, $t\in[\TNULL,\TEINS]$ can be regarded as an element
  of $\widehat{W}_{\widehatNeumann}^{-1,q}$. We denote $\Rbvarphi$
  as a function from $[\TNULL,\TEINS]$ into
  $\widehat{W}_{\widehatNeumann}^{-1,q}$ by $\rbvarphi$. 
  The H\"older continuous differentiability of $\Rbvarphi$ entails 
  the H\"older continuous differentiability of
  \begin{math}
    \rbvarphi
    \from{[\TNULL,\TEINS] \to \widehat{W}_{\widehatNeumann}^{-1,q}}
  \end{math}
  with the same H\"older exponent. 
\end{rem}

\subsection{The linear Poisson equation}

Let us assume the following about $\doping$ --- the doping
profile (or control parameter) on the right hand side of Poisson's
equation \eqref{eq:poi}.
\begin{assu} 
  \label{assu:dop}
  The function 
  \begin{math}
    \doping
    \from{[\TNULL,\TEINS] \to \widehat{W}_{\widehatNeumann}^{-1,q}}
  \end{math}
  is continuously differentiable with H\"older continuous derivative. 
  We define a ``generalized doping'' 
  \begin{equation} 
    \label{e-dop}
    \immotile\from{[\TNULL,\TEINS] \to \widehat{W}_{\widehatNeumann}^{-1,q}}
    \qquad
    \text{by}
    \qquad
    \immotile(t) 
    \df 
    \doping(t) 
    + 
    \rbvarphi(t)
    ,
    \quad
    t\in[\TNULL,\TEINS]
    .
  \end{equation}
\end{assu}
We now define what is a solution of Poisson's equation \eqref{eq:poi}.
\begin{defn}
  \label{def:linpoisson}
  Let $u_k\in\widehat{W}_{\widehatNeumann}^{-1,q}$, $k=1,2$ be given. 
  We say that $\Varphi$ is a solution of Poisson's equation
  \eqref{eq:poi} at $t\in[\TNULL,\TEINS]$, if
  \begin{equation} 
    \label{eq:splitoff}
    \Varphi = \varphi + \dbvarphi(t),
  \end{equation}
  and $\varphi\in\widehat{W}_{\widehatNeumann}^{1,q}$ is the unique
  solution of
  \begin{equation} 
    \label{e-linpoi}
    \Linpoisson\varphi = \immotile(t) + u_1 - u_2.
  \end{equation}
  $\varphi$ and $\Varphi$ depend parametrically on $t$, $u_1$, and
  $u_2$. If convenient, we indicate the dependence on $t$ by writing
  $\varphi(t)$ and $\Varphi(t)$, respectively.
\end{defn}
\begin{rem} 
  \label{rem:boundjustify}
  With respect to the boundary conditions in \eqref{eq:poi} it should
  be noted that \eqref{eq:trace2} and the property
  $\varphi\in\widehat{W}_{\widehatNeumann}^{1,q}$ give
  \begin{math}
    \Varphi|_{\widehatDirichlet} = \Dbvarphi.
  \end{math}
  Additionally, if $\doping$, $u_1$, and $u_2$ belong to the space
  $\widehat{L}^1$, then \eqref{e-dop}, \eqref{eq:splitoff} and
  \eqref{e-linpoi} together with \eqref{eq:trace1} imply
\begin{math}
  {\nu}\cdot{ 
    \left(
      \varepsilon \grad \Varphi
    \right)
  }
  + 
  \varepsilon_{\widehatNeumann} \Varphi 
  =\Rbvarphi(t)
  ,
\end{math} 
\cf\ \eg\ 
\cite[Ch.~1.2]{ciarlet79} or
\cite[Ch.~II.2]{gajewski:groeger:etal74}.
\end{rem}
Throughout this section we demand several times H\"older continuity of
functions and/or their derivatives. Clearly, there is a common
H\"older exponent which we will denote from now on by $\eta$.

\section[Precise Formulation of the Problem]
  {Precise Formulation of the Problem}
\label{sec:solution}
We are now going to define the problem outlined in
\secref{sec:setting}.
\begin{defn} 
  \label{def:vanroos}
  We say the van Roosbroeck system admits a local in time solution, if
  there is a time $T\in]\TNULL,\TEINS]$ and 
  \begin{math}
    (\Varphi, \Fullphi )=(\Varphi, {\Fullphi}_1,{\Fullphi}_2) 
  \end{math}
  such that
  \begin{equation}
    \label{eq:initial-condition}
    \Fullphi(\TNULL) 
    = (\Fullphi_1(\TNULL),\Fullphi_2(\TNULL))
    =(\iniphi_1,\iniphi_2) 
    \in \xoplus{W}^{1,q},
  \end{equation}
  \begin{equation} 
    \label{eq:sol-varphi}
    \varphi \df \Varphi-\dbvarphi 
    \in 
    C([\TNULL,T];\widehat{W}^{1,q}_{\widehatNeumann}) 
    \cap 
    C^1(]\TNULL,T[;\widehat{W}^{1,q}_{\widehatNeumann})
  \end{equation}
  \begin{equation}
    \label{eq:sol-phi}
    \phi \df    \Fullphi - \bvphi{}
    \in
    C^1(]\TNULL,T[,\xoplus{L}^p) 
    \cap 
    C(]\TNULL,T],\domainA)
    \cap 
    C([\TNULL,T],[\xoplus{L}^p,\domainA]_\theta),
  \end{equation}
  fulfill the Poisson equation and the current continuity equations: 
  \begin{equation}
    \label{eq:poisson}
    \Linpoisson(\varphi(t)) 
    = \immotile(t) + \extension u_1(t) - \extension u_2(t)
    \quad
    t\in[\TNULL,T]
    ,
  \end{equation}
  \begin{equation}
    \label{eq:continuity}
    u_k'(t) - \dive j_k(t) 
    = 
    \reaction_k(t,\Varphi(t),\Fullphi(t)),
    \quad k=1,2,
    \quad t \in ]\TNULL,T[
    .
  \end{equation}   
  The carrier densities and the current densities are given by
  \begin{align}
    \label{eq:car-density}
    u_k(t) 
    &\df 
    \rho_k(t) \mathcal{F}_k
    \big(
    \chempot_k(t) 
    \big),
    \\
    \label{eq:cur-density}
    j_k(t) 
    &\df 
    \mathcal{G}_k\big(
    \chempot_k(t) 
    \big)
    \mu_k\grad {\Fullphi}_k(t),
    \\
    \label{eq:chemical-potential}
    \chempot_k(t) 
    &\df 
    {\Fullphi}_k(t) + (-1)^k\restrict \Varphi(t) + b_k(t).
  \end{align}
  and satisfy
  \begin{equation}
    \label{eq:ukjk-regularity}
    u_k \in C([\TNULL,T],L^\infty) \cap C^1(]\TNULL,T[,L^p),
  \end{equation}
  \begin{equation}
    \label{eq:jk-regularity}
    j_k \in C([\TNULL,T],{L}^q),
  \end{equation}
  \begin{equation}
    \label{eq:divjk-regularity}
    \dive j_k \in C(]\TNULL,T],{L}^p)
  \end{equation}
  for $k=1,2$.     
\end{defn}

\section[Reformulation as a quasi-linear parabolic system]
{Reformulation as a quasi-linear parabolic system}
\label{sec:nl-reform}

In this section we provide the tools to rewrite
the problem from \defnref{def:vanroos}
%%Problem~\eqref{eq:initial-condition}--\eqref{eq:chemical-potential} 
as a quasi-linear system for the continuity equations. To that end we
eliminate the electrostatic potential from the continuity equations.
Replacing the carrier densities $u_1$ and $u_2$ on the right hand side
of \eqref{eq:poisson} by \eqref{eq:car-density} making use of
\eqref{eq:chemical-potential} and \eqref{eq:splitoff} one obtains a
nonlinear Poisson equation for $\varphi$. We solve this equation with
respect to prescribed parameters $b_k$ and ${\Fullphi}_k$, $k=1,2$,
which we will assume here to be from $L^\infty$. This way to decouple
van Roosbroeck's equations into a nonlinear Poisson equation and a
system of parabolic equations is also one of the fundamental
approaches to the numerical solution of the van Roosbroeck system. It
is due to Gummel \cite{gummel64} and was the first reliable numerical
technique to solve these equations for carriers in an operating
semiconductor device structure.

\subsection{The nonlinear Poisson equation}
\label{sec:nl-poisson}

We are now going to prove the unique solvability of the nonlinear
Poisson equation and some properties of its solution.  First we show
that the supposed admissibility of the carrier distribution functions
$\mathcal{F}_k$ ensures that the relation between a potential and its
corresponding carrier density is monotone and even continuously
differentiable when considered between adequate spaces.
\begin{lem} 
  \label{lem:cardens}
  Let $\rho$ and $g$ be from ${L}^\infty$ and
  $\mathcal{F}=\mathcal{F}_k$ be an admissible carrier distribution
  function, see \assuref{assu:distri}.
  \par
  1.~The operator
  \begin{equation} 
    \label{eq:cardens}
    \widehat{W}_{\widehatNeumann}^{1,2} 
    \ni 
    h
    \longmapsto 
    \extension \rho \mathcal{F}(g +\restrict h) 
    \in \widehat{L}^2
  \end{equation}
  is well defined, continuous and bounded. Its composition with the
  embedding
  \begin{math}
    \widehat{L}^2 \embedcontinuously \widehat{W}_{\widehatNeumann}^{-1,2}
  \end{math}
  is monotone.
  \par
  2.~The Nemyckii operator 
  \begin{equation*} 
    L^\infty 
    \ni 
    h
    \longmapsto 
    \rho {\mathcal{F}}(g + \restrict h) 
  \end{equation*}
  induced by the function
  \begin{equation*}
    \Omega \times \mathbb{R} 
    \ni 
    (x,s) 
    \longmapsto 
    \rho(x) \mathcal{F}(g(x)+s),
  \end{equation*}
  maps $L^\infty$ continuously into itself and is even continuously
  differentiable.  Its Fr{\'e}chet derivative at $h\in{L}^\infty$ is
  the multiplication operator given by the essentially bounded
  function
  \begin{equation}
    \label{eq:cardens:deriv}
    \Omega \ni x \longmapsto  \rho(x) \mathcal{F}^\prime(g(x)+ h(x)).
  \end{equation}
\end{lem}
\begin{proof} 
  Indeed, the assumption that the carrier distribution functions
  should be admissible assures that the operator \eqref{eq:cardens} is
  well defined, continuous and bounded, \cf\ \cite{trudinger:67} for
  the case of an exponential, and \cf\ \cite[Chapter~3]{appell90} for
  the case of a polynomially bounded function.  The asserted
  monotonicity follows from the monotonicity of the function
  $\mathcal{F}$ and the fact that the duality between
  $\widehat{W}_{\widehatNeumann}^{1,2}$ and
  $\widehat{W}_{\widehatNeumann}^{-1,2}$ is the extension of the
  $\widehat{L}^2$ duality:
  \begin{multline*}
    \dual{
      \extension  \rho \mathcal{F} ( g + \restrict h_1)
      -
      \extension \rho \mathcal{F} ( g + \restrict h_2 )
    }{h_1 - h_2}_{\widehat{W}_{\widehatNeumann}^{1,2}}
    \\
    =
    \int_{\widehat{\Omega}}
    \left(
      \extension \rho \mathcal{F} ( g + \restrict h_1 )
      -
      \extension \rho \mathcal{F} ( g + \restrict h_2 )
    \right)
    \left(
      h_1 - h_2
    \right)
    \der{x}\\
    = 
    \int_{\Omega}
    \left(
      \rho \mathcal{F} ( g + \restrict h_1 )
      -
      \rho \mathcal{F} ( g + \restrict h_2 )
    \right)
    \left(
      \restrict  h_1 - \restrict h_2
    \right)
    \der{x}
    \ge 0
    \;\;
    \text{for all $h_1$, $h_2 \in \widehat{W}_{\widehatNeumann}^{1,2}$.}
  \end{multline*}
  \par
  The second assertion follows from a result by Gr\"oger and Recke, \cf\
  \cite[Thm~5.1]{recke-groeger-nodea}.
\end{proof}
\begin{cor} 
  \label{cor:implica}
  The mapping
  \begin{equation*} 
    \widehat{W}^{1,q} 
    \ni 
    h
    \longmapsto 
    \extension \rho {\mathcal{F}}(g + \restrict h) 
  \end{equation*}
  takes its values in $\widehat{L}^\infty$ and is also continuously
  differentiable. Its derivative at a point $h\in\widehat{W}^{1,q}$
  equals the multiplication operator which is induced by the function
  $\extension\rho\mathcal{F}'(g+\restrict h)$.
\end{cor}
\begin{thm} 
  \label{thm:monotone}
  Under \assuref{assu:distri} on the distribution functions
  $\mathcal{F}_1$, $\mathcal{F}_2$ and \assuref{assu:effband} the
  following statements are true:
  \par
  1.~For any pair of functions 
  $z=(z_1,z_2)\in\xoplus{L}^\infty$
  the operator 
  \begin{equation} 
    \label{eq:monmap}
    \varphi 
    \longmapsto 
    \Linpoisson \varphi
    - \extension \rho_1 \mathcal{F}_1 (z_1 - \restrict \varphi)
    + \extension \rho_2 \mathcal{F}_2 (z_2 + \restrict \varphi)
  \end{equation}
  is strongly monotone and continuous from
  $\widehat{W}_{\widehatNeumann}^{1,2}$ to
  $\widehat{W}_{\widehatNeumann}^{-1,2}$, where the operator
  $\Linpoisson$ is according to \defnref{d-opera}. The monotonicity
  constant of \eqref{eq:monmap} is a least that of $\Linpoisson$.
  \par
  2.~For all
  $f\in\widehat{W}_{\widehatNeumann}^{-1,2}$ and
  $z=(z_1,z_2)\in\xoplus{L}^\infty$ 
  the nonlinear Poisson equation
  \begin{equation}
    \label{eq:nlp}
    \Linpoisson \varphi 
    - \extension \rho_1 \mathcal{F}_1(z_1-\restrict \varphi)
    + \extension \rho_2 \mathcal{F}_2(z_2+\restrict \varphi)
    =f
  \end{equation}
  admits exactly one solution $\varphi$ which we denote by
  $\nlpsolution(f,z)$.  This solution belongs to
  $\widehat{W}_{\widehatNeumann}^{1,2}$ and satisfies the estimate
  \begin{equation*}
    \norm{\varphi}_{\widehat{W}_{\widehatNeumann}^{1,2}} 
    \le 
    \frac{1}{m} 
    \lrnorm{
      \extension \rho_1 \mathcal{F}_1(z_1)
      - 
      \extension \rho_2 \mathcal{F}_2(z_2) 
      +
      f
    }_{\widehat{W}_{\widehatNeumann}^{-1,2}},
  \end{equation*}
  where $m$ is the monotonicity constant of $\Linpoisson$.
  \par 
  3.~The maximal restriction of the operator \eqref{eq:monmap} to the
  range space $\widehat{W}_{\widehatNeumann}^{-1,q}$ has the domain
  $\widehat{W}_{\widehatNeumann}^{1,q}$. Moreover, if
  $M$ is a bounded subset of
  \begin{math}
    \widehat{W}_{\widehatNeumann}^{-1,q}
    \oplus \xoplus{L}^\infty,
  \end{math}
  then the set 
  \begin{math}
    \left \{
      \nlpsolution(f,z) 
      \with (f, z) \in M 
    \right \}
  \end{math}
  is bounded in $\widehat{W}_{\widehatNeumann}^{1,q}$.
  \par
  4.~The mapping 
  \begin{math}
    \nlpsolution
    \from{
      \widehat{W}_{\widehatNeumann}^{-1,q} \oplus \xoplus{L}^\infty
      \to
      \widehat{W}_{\widehatNeumann}^{1,q}
    }
  \end{math}
  is continuously differentiable. Let $(F,Z)=(F,Z_1,Z_2)$ be from 
  \begin{math}
    \widehat{W}_{\widehatNeumann}^{-1,q} \oplus \xoplus{L}^\infty
    ;
  \end{math}
  we define the function
  \begin{equation} 
    \label{eq:pk}
    \Lindensity_k 
    \df 
    \extension \rho_k 
    \mathcal{F}'_k (Z_k  + (-1)^k\restrict \nlpsolution(F,Z))
    ,
  \end{equation}
  and we also denote the corresponding multiplication operator on
  $\widehat{\Omega}$ by $\Lindensity_k$.
  Then the Fr{\'e}chet derivative
  \begin{math}
    \partial \nlpsolution 
  \end{math}
  at a point 
  \begin{math}
    (F,Z) = (F,Z_1,Z_2)
  \end{math}
  is the bounded linear mapping given by
  \begin{equation} 
    \label{eq:nlp-fderiv}
    \left[
      \partial \nlpsolution(F,Z)
    \right]
    (f,z)
    = 
    \left(
      \Linpoisson + \Lindensity_1 + \Lindensity_2
    \right)^{-1}
    \left(
      f 
      + \Lindensity_1  \extension z_1 
      - \Lindensity_2  \extension z_2 
    \right), 
    \quad k=1,2
  \end{equation}
  for all
  \begin{math}
    (f,z) = (f,(z_1,z_2))
    \in 
    \widehat{W}_{\widehatNeumann}^{-1,q} 
    \oplus 
    \xoplus{L}^\infty
    .
  \end{math}
  \par
  5.~The norm of 
  \begin{math}
    \partial \nlpsolution(F,Z) 
    \in  
    \mathcal{B}(
    \widehat{W}_{\widehatNeumann}^{-1,q} \oplus \xoplus{L}^\infty
    ;
    \widehat{W}_{\widehatNeumann}^{1,q}
    )
  \end{math}
  can be estimated as follows:
  \begin{multline*}
    \norm{\partial \nlpsolution(F,Z)}_{
      \mathcal{B}(
      \widehat{W}_{\widehatNeumann}^{-1,q} \oplus \xoplus{L}^\infty
      ;
      \widehat{W}_{\widehatNeumann}^{1,q}
      )}
    \\
    \le
    2
    \norm{\Linpoisson^{-1}}_{
      \mathcal{B}(L^2;\widehat{W}_{\widehatNeumann}^{1,q})}
    \sqrt{
      \norm{\Lindensity_1 + \Lindensity_2}_{L^\infty}
      \norm{\Lindensity_1 + \Lindensity_2}_{L^1}
    }
    +
    \norm{\Linpoisson^{-1}}_{\mathcal{B}(
      \widehat{W}_{\widehatNeumann}^{-1,q};
      \widehat{W}_{\widehatNeumann}^{1,q})}
    \\
    +
    \norm{\Linpoisson^{-1}}_{\mathcal{B}(\widehat{L}^2;
      \widehat{W}_{\widehatNeumann}^{1,q})}
    \sqrt{
      \norm{\Lindensity_1 + \Lindensity_2}_{L^\infty}
    }
    \norm{\Linpoisson^{-1/2}}_{\mathcal{B}(
      \widehat{W}_{\widehatNeumann}^{-1,q};\widehat{L}^2)}
  \end{multline*}
\end{thm}
\begin{proof}
  1.~The assumption that $\widehatDirichlet$ is not empty or
  $\varepsilon_{\widehatNeumann}$ is positive on a set of positive arc
  measure ensures that the operator $\Linpoisson$ is strongly
  monotone.  Thus, taking into account \lemref{lem:cardens}, the
  mapping \eqref{eq:monmap} is strongly monotone and continuous from
  $\widehat{W}_{\widehatNeumann}^{1,2}$ to
  $\widehat{W}_{\widehatNeumann}^{-1,2}$.
  
  2.~The second assertion follows from the first one by standard
  results on monotone operators, \cf\ \eg\ 
  \cite{gajewski:groeger:etal74}.
  
  3.~For $f\in\widehat{W}^{-1,2}_{\widehatNeumann}$ the solution
  $\nlpsolution(f,z)$ is from $\widehat{W}^{1,2}_{\widehatNeumann}$
  and hence,
  \begin{equation*}
    - \extension \rho_1 
    \mathcal{F}_1\big(z_1- \restrict \nlpsolution(f,z)\big)
    + \extension \rho_2 
    \mathcal{F}_2\big(z_2+ \restrict \nlpsolution(f,z)\big)
    \in
    \widehat{L}^2
    \embedcontinuously
    \widehat{W}^{-1,q}_{\widehatNeumann},
  \end{equation*}
  \cf\ \lemref{lem:cardens}.  By the second assertion of the theorem,
  the set 
  \begin{equation*}
    \left\{\nlpsolution(f,z){\with}(f,z){\in}M\right\}
    \quad
    \text{is bounded in $\widehat{W}^{1,2}_{\widehatNeumann}$.}
  \end{equation*}
  From this we conclude again by \lemref{lem:cardens} that the set
  \begin{equation*}
    \left\{
      \extension \rho_1 
      \mathcal{F}_1\big(z_1-\restrict \nlpsolution(f,z)\big)
      -
      \extension \rho_2 
      \mathcal{F}_2\big(z_2+\restrict  \nlpsolution(f,z)\big)
      \with
      (f,z) \in M
    \right\}
  \end{equation*}
  is bounded in $\widehat{L}^2$, and hence, is bounded in
  $\widehat{W}^{-1,q}_{\widehatNeumann}$. Thus, the set
  \begin{equation*}
    \left\{
      \extension \rho_1 
      \mathcal{F}_1\big(z_1-\restrict \nlpsolution(f,z)\big)
      -
      \extension \rho_2 
      \mathcal{F}_2\big(z_2+\restrict  \nlpsolution(f,z)\big)
      +
      f
      \with
      (f,z) \in M
    \right\}
  \end{equation*}
  is also bounded in $\widehat{W}^{-1,q}_{\widehatNeumann}$.
  Consequently, the image of this set under $\Linpoisson^{-1}$ is
  bounded in $\widehat{W}^{1,q}_{\widehatNeumann}$.
  
  4.~We define an auxiliary mapping
  \begin{math}
    \mathcal{K}\from{
      \widehat{W}_{\widehatNeumann}^{1,q}
      \oplus
      \widehat{W}_{\widehatNeumann}^{-1,q}
      \oplus
      \xoplus{L}^\infty
      \to
      \widehat{W}_{\widehatNeumann}^{-1,q}
    }
  \end{math}
  by
  \begin{equation*} 
    \mathcal{K}(\varphi,f,z) 
    \df
    \Linpoisson \varphi
    - \extension \rho_1
    \mathcal{F}_1(z_1 - \restrict \varphi)
    + \extension \rho_2
    \mathcal{F}_2(z_2 + \restrict \varphi)
    - f
  \end{equation*}
  such that
  \begin{math}
    \mathcal{K}\big(\nlpsolution(f,z),f,z\big) 
    =
    0
  \end{math}
  for all $f\in\widehat{W}_{\widehatNeumann}^{-1,q}$ and all
  $z\in\xoplus{L}^\infty$.  The assertion follows from the
  \thmtitle{Implicit Function Theorem} if we can prove that
  $\mathcal{K}$ is continuously differentiable and the partial
  derivative with respect to $\varphi$ is a topological isomorphism
  between $\widehat{W}_{\widehatNeumann}^{1,q}$ and
  $\widehat{W}_{\widehatNeumann}^{-1,q}$.  
  For any
  \begin{math}
    \varphi \in \widehat{W}_{\widehatNeumann}^{1,q},
  \end{math}
  \begin{math}
    f \in \widehat{W}_{\widehatNeumann}^{-1,q},
  \end{math}
  and
  \begin{math}
    z \in \xoplus{L}^\infty
  \end{math}
  the partial derivatives of $\mathcal{K}$ are given by
  \begin{eqnarray}
    \label{eq:grad-K-varphi}
    \partial_\varphi
    \mathcal{K}(\varphi,f,z)
    & = &
    \Linpoisson
    +
    \sum_{k=1}^2
    \extension \rho_k \mathcal{F}'_k
    (z_k 
    + (-1)^k
    \restrict \varphi)
    \in
    \mathcal{B}(
    \widehat{W}_{\widehatNeumann}^{1,q};
    \widehat{W}_{\widehatNeumann}^{-1,q}
    ),
    \\
    \label{eq:grad-K-nlpright}
    \partial_f
    \mathcal{K}(\varphi,f,z)
    & = &
    -\id 
    \in
    \mathcal{B}(
    \widehat{W}_{\widehatNeumann}^{-1,q};
    \widehat{W}_{\widehatNeumann}^{-1,q}
    ),
    \\
    \label{eq:grad-K-nlpfixed}
    \partial_{z_k}
    \mathcal{K}(\varphi,f,z)
    & = &
    (-1)^k
    \extension \rho_k \mathcal{F}'_k
    (z_k 
    + (-1)^k
    \restrict \varphi)
    \in
    \widehat{L}^\infty
    \embedcontinuously
    \mathcal{B}(
    \xoplus{L}^\infty;
    \widehat{W}_{\widehatNeumann}^{-1,q}
    )
  \end{eqnarray}
  and they are continuous, \cf\ \lemref{lem:cardens} and
  \cite[\S5]{recke-groeger-nodea}.
 
  Now we consider the equation
  \begin{equation} 
    \label{eq:laxmil}
    \Linpoisson \psi
    +
    \sum_{k=1}^2
    \extension \rho_k 
    \mathcal{F}'_k
    (z_k + (-1)^k \restrict \varphi)
    \psi
    =
    f \in \widehat{W}_{\widehatNeumann}^{-1,q}
  \end{equation}                               %
  Because
  \begin{math}
    \sum_{k=1}^2
    \extension \rho_k 
    \mathcal{F}'_k
    (z_k + (-1)^k \restrict \varphi)
  \end{math}
  is a positive function from $\widehat{L}^\infty$, \eqref{eq:laxmil}
  has exactly one solution
  $\psi\in\widehat{W}_{\widehatNeumann}^{1,2}$ by the
  \thmtitle{Lax-Milgram-Lemma}. Moreover,  
  \begin{equation*}
    \sum_{k=1}^2
    \extension \rho_k 
    \mathcal{F}'_k
    (z_k + (-1)^k \restrict \varphi)
    \psi
    \in
    \widehat{L}^2 
    \embedcontinuously
    \widehat{W}_{\widehatNeumann}^{-1,q}
    ,
  \end{equation*}
  and
  \begin{math}
    \Linpoisson\from{
      \widehat{W}_{\widehatNeumann}^{1,q}
      \to
      \widehat{W}_{\widehatNeumann}^{-1,q}
    }
  \end{math}
  is a topological isomorphism. Thus, a rearrangement of terms in
  \eqref{eq:laxmil} gives
  $\psi\in\widehat{W}_{\widehatNeumann}^{1,q}$.

  5.~We now estimate the Fr{\'e}chet derivative \eqref{eq:nlp-fderiv}:
  \begin{multline}
    \label{eq:impl02}
    \lrnorm{
      (\Linpoisson
      +\Lindensity_1
      +\Lindensity_2)^{-1}
      (f
      +\Lindensity_1 \extension z_1
      -\Lindensity_2  \extension z_2
      )
    }_{\widehat{W}_{\widehatNeumann}^{1,q}}
    \\
    \le
    \lrnorm{
      (\Linpoisson
      +\Lindensity_1
      +\Lindensity_2)^{-1}
      f
    }_{\widehat{W}_{\widehatNeumann}^{1,q}}
    \\
    +
    \lrnorm{
      (\Linpoisson
      +\Lindensity_1
      +\Lindensity_2)^{-1}
      (\Lindensity_1 \extension z_1
      -\Lindensity_2 \extension z_2
      )
    }_{\widehat{W}_{\widehatNeumann}^{1,q}}
    .
  \end{multline}
  We treat the right hand side terms separately;
  for the second addend one obtains
  \begin{multline}
    \label{eq:impl02aa}
    \lrnorm{
      (\Linpoisson
      +\Lindensity_1
      +\Lindensity_2)^{-1}
      (\Lindensity_1 \extension z_1
      -\Lindensity_2 \extension z_2
      )
    }_{\widehat{W}_{\widehatNeumann}^{1,q}}
    \\
    \le
    \lrnorm{
      (\Linpoisson
      +\Lindensity_1
      +\Lindensity_2)^{-1}
      \sqrt{\Lindensity_1 + \Lindensity_2}
    }_{\mathcal{B}(\widehat{L}^2;\widehat{W}_{\widehatNeumann}^{1,q})}
    \lrnorm{g}_{L^2}
    ,
  \end{multline}
  where the function $g\in{L^2}$ is defined by 
  \begin{equation}
    \label{eq:g}
    g(x) \df
    \frac{\Lindensity_1(x)z_1(x)
      -\Lindensity_2(x)z_2(x)}
    {\sqrt{\Lindensity_1(x)+\Lindensity_2(x)}}
    \quad\text{for $x\in\Omega$.}
  \end{equation}
%%   \begin{equation*}
%%     g(x) \df
%%     \begin{cases}
%%       0 
%%       & \text{if $\Lindensity_1(x)=\Lindensity_2(x)=0$,}
%%       \\
%%       \frac{\Lindensity_1(x)z_1(x)
%%         -\Lindensity_2(x)z_2(x)}
%%       {\sqrt{\Lindensity_1(x)+\Lindensity_2(x)}}
%%       & \text{elsewhere on $\Omega$}
%%     \end{cases}
%%   \end{equation*}
  Please note that the functions $\Lindensity_k$ are strictly positive
  almost everywhere in $\Omega$ due to the positivity of the
  distribution functions and \assuref{assu:effband}.  For the function
  $g$ in \eqref{eq:g} one has the following bound:
  \begin{equation*}
    \norm{g}_{L^2}
    \le 
    \sqrt{\norm{\Lindensity_1+\Lindensity_2}_{\widehat{L}^1}}
    \left(
      \norm{z_1}_{L^\infty}
      +
      \norm{z_2}_{L^\infty} 
    \right)
    .
  \end{equation*}
  Making use of the operator identity
  \begin{equation}
    \label{eq:opid}
    (\Linpoisson+\Lindensity_1+\Lindensity_2)^{-1}
    = 
    \Linpoisson^{-1}
    -
    \Linpoisson^{-1}
    (\Lindensity_1+\Lindensity_2)
    (\Linpoisson+\Lindensity_1+\Lindensity_2)^{-1}
  \end{equation}
  one obtains
  \begin{multline*}
    \lrnorm{
      (\Linpoisson
      +\Lindensity_1
      +\Lindensity_2)^{-1}
      \sqrt{\Lindensity_1 + \Lindensity_2}
    }_{\mathcal{B}(\widehat{L}^2;\widehat{W}_{\widehatNeumann}^{1,q})}
    \le
    \lrnorm{
      \Linpoisson^{-1}
      \sqrt{\Lindensity_1 + \Lindensity_2}
    }_{\mathcal{B}(\widehat{L}^2;\widehat{W}_{\widehatNeumann}^{1,q})}
    \\
    +
    \lrnorm{    
      \Linpoisson^{-1}
      \sqrt{\Lindensity_1 + \Lindensity_2}
      \sqrt{\Lindensity_1 + \Lindensity_2}
      (\Linpoisson
      +\Lindensity_1
      +\Lindensity_2)^{-1}
      \sqrt{\Lindensity_1 + \Lindensity_2}
    }_{\mathcal{B}(\widehat{L}^2;\widehat{W}_{\widehatNeumann}^{1,q})}
    \\
    \le
    \lrnorm{
      \Linpoisson^{-1}
    }_{\mathcal{B}(\widehat{L}^2;\widehat{W}_{\widehatNeumann}^{1,q})}
    \sqrt{\lrnorm{\Lindensity_1 + \Lindensity_2}_{\widehat{L}^\infty}}
    \;\times
    \\
    \times
    \left(
      1
      +
      \lrnorm{    
        \sqrt{\Lindensity_1 + \Lindensity_2}
        (\Linpoisson
        +\Lindensity_1
        +\Lindensity_2)^{-1/2}
      }^2_{\mathcal{B}(\widehat{L}^2)}
    \right)
  \end{multline*}
  We note that
  \begin{equation} 
    \label{eq:subord}
    \lrnorm{    
      \sqrt{\Lindensity_1 + \Lindensity_2}
      (\Linpoisson
      +\Lindensity_1
      +\Lindensity_2)^{-1/2}
    }_{\mathcal{B}(\widehat{L}^2)}
    \le 1
  \end{equation}
  because the bounded multiplication operator $\Lindensity_1+\Lindensity_2$
  is form subordinated to 
  \begin{math}
    \Linpoisson
    +
    \Lindensity_1 + \Lindensity_2,
  \end{math}
  \cf\ \eg\ \cite[VI.2.6]{kato}.
  Thus, we get for the second addend of \eqref{eq:impl02}:
  \begin{multline}
    \label{eq:impl02a}
    \lrnorm{
      (\Linpoisson
      +\Lindensity_1
      +\Lindensity_2)^{-1}
      (\Lindensity_1\extension z_1
      -\Lindensity_2\extension z_2
      )
    }_{\widehat{W}_{\widehatNeumann}^{1,q}}
    \\
    \le
    2\,
    \lrnorm{
      \Linpoisson^{-1}
    }_{\mathcal{B}(\widehat{L}^2;\widehat{W}_{\widehatNeumann}^{1,q})}
    \sqrt{\lrnorm{\Lindensity_1 + \Lindensity_2}_{\widehat{L}^\infty}}
    \sqrt{\lrnorm{\Lindensity_1 + \Lindensity_2}_{\widehat{L}^1}}
    \left(
      \norm{z_1}_{L^\infty}
      +
      \norm{z_2}_{L^\infty} 
    \right)
  \end{multline}
  Applying \eqref{eq:opid} to the first term on the right hand side of
  \eqref{eq:impl02} we find
  \begin{multline}
    \label{eq:impl02bb}
    \lrnorm{
      (\Linpoisson
      +\Lindensity_1
      +\Lindensity_2)^{-1}
      f
    }_{\widehat{W}_{\widehatNeumann}^{1,q}}
    \le
    \lrnorm{
      \Linpoisson^{-1}
    }_{\mathcal{B}(
      \widehat{W}_{\widehatNeumann}^{-1,q}
      ;
      \widehat{W}_{\widehatNeumann}^{1,q})}
    \lrnorm{
      f
    }_{\widehat{W}_{\widehatNeumann}^{-1,q}}
    \\
    +
    \lrnorm{
      \Linpoisson^{-1}
    }_{\mathcal{B}(\widehat{L}^2;\widehat{W}_{\widehatNeumann}^{1,q})}
    \lrnorm{
      (\Lindensity_1+\Lindensity_2)
      (\Linpoisson+\Lindensity_1+\Lindensity_2)^{-1}
    }_{\mathcal{B}(\widehat{W}_{\widehatNeumann}^{-1,q};\widehat{L}^2)}
    \lrnorm{
      f
    }_{\widehat{W}_{\widehatNeumann}^{-1,q}}.
  \end{multline}
  The terms
  \begin{math}
    \lrnorm{
      \Linpoisson^{-1}
    }_{\mathcal{B}(
      \widehat{W}_{\widehatNeumann}^{-1,q}
      ;
      \widehat{W}_{\widehatNeumann}^{1,q})}
  \end{math}
  and
  \begin{math}
    \lrnorm{
      \Linpoisson^{-1}
    }_{\mathcal{B}(\widehat{L}^2;\widehat{W}_{\widehatNeumann}^{1,q})}
  \end{math}
  are finite.
  As for the remaining term
  \begin{multline*}
    \lrnorm{
      (\Lindensity_1+\Lindensity_2)
      (\Linpoisson+\Lindensity_1+\Lindensity_2)^{-1}
    }_{\mathcal{B}(\widehat{W}_{\widehatNeumann}^{-1,q};\widehat{L}^2)}
    \\
    \le
    \sqrt{\lrnorm{\Lindensity_1 + \Lindensity_2}_{\widehat{L}^\infty}}
    \lrnorm{    
      \sqrt{\Lindensity_1 + \Lindensity_2}
      (\Linpoisson
      +\Lindensity_1
      +\Lindensity_2)^{-1/2}
    }_{\mathcal{B}(\widehat{L}^2)}
    \\
    \lrnorm{    
      (\Linpoisson
      +\Lindensity_1
      +\Lindensity_2)^{-1/2}
      \Linpoisson^{1/2}
    }_{\mathcal{B}(\widehat{L}^2)}
    \lrnorm{  
      \Linpoisson^{-1/2} 
    }_{\mathcal{B}(\widehat{W}_{\widehatNeumann}^{-1,q};\widehat{L}^2)}
  \end{multline*}
  we note that 
  \begin{math}
    \lrnorm{
      \Linpoisson^{-1/2}
    }_{\mathcal{B}(\widehat{W}_{\widehatNeumann}^{-1,q};\widehat{L}^2)}
  \end{math}
  is finite, 
  since 
  \begin{math}
    \widehat{W}_{\widehatNeumann}^{-1,q} 
  \end{math}
  embeds continuously into
  \begin{math}
    \widehat{W}_{\widehatNeumann}^{-1,2}
  \end{math}
  and
  \begin{math}
    \Linpoisson^{1/2}
    \from{\widehat{L}^2 \to \widehat{W}_{\widehatNeumann}^{-1,2}}
  \end{math}
  is a topological isomorphism.
                                %  \cf\ \defnref{def:poilin}.
  Again, $\Linpoisson$ is form subordinated to
  $\Linpoisson+\Lindensity_1+\Lindensity_2$.  Hence, besides
  \eqref{eq:subord} one has
  \begin{equation*}
%   \norm{    
%      \sqrt{\Lindensity_1 + \Lindensity_2}
%      (\Linpoisson 
%      +\Lindensity_1
%      +\Lindensity_2)^{-1/2}
%    }_{\mathcal{B}(\widehat{L}^2)}
%    \le 1
%    ,\qquad
    \norm{    
      (\Linpoisson
      +\Lindensity_1
      +\Lindensity_2)^{-1/2}
      \Linpoisson^{1/2}
    }_{\mathcal{B}(\widehat{L}^2)}
    \le 1
    .
  \end{equation*}
  Thus, we get from \eqref{eq:impl02bb}:
  \begin{multline}
    \label{eq:impl02b}
    \lrnorm{
      (\Linpoisson
      +\Lindensity_1
      +\Lindensity_2)^{-1}
      f
    }_{\widehat{W}_{\widehatNeumann}^{1,q}}
    \le
    \lrnorm{
      \Linpoisson^{-1}
    }_{\mathcal{B}(
      \widehat{W}_{\widehatNeumann}^{-1,q}
      ;
      \widehat{W}_{\widehatNeumann}^{1,q})}
    \lrnorm{
      f
    }_{\widehat{W}_{\widehatNeumann}^{-1,q}}
    \\
    +
    \lrnorm{
      \Linpoisson^{-1}
    }_{\mathcal{B}(\widehat{L}^2;\widehat{W}_{\widehatNeumann}^{1,q})}
    \sqrt{\lrnorm{\Lindensity_1 + \Lindensity_2}_{\widehat{L}^\infty}}
    \lrnorm{  
      \Linpoisson^{-1/2} 
    }_{\mathcal{B}(\widehat{W}_{\widehatNeumann}^{-1,q};\widehat{L}^2)}
    \lrnorm{
      f
    }_{\widehat{W}_{\widehatNeumann}^{-1,q}}.
  \end{multline}
  Inserting \eqref{eq:impl02a} and \eqref{eq:impl02b} into
  \eqref{eq:impl02} finishes the proof.
\end{proof}
\begin{cor} 
  \label{cor:boundedlip}
  Let the assumptions of \thmref{thm:monotone} be satisfied. Then
  holds true:
  
  1.~The mapping 
  \begin{math}
    \nlpsolution
    \from{
      \widehat{W}_{\widehatNeumann}^{-1,q} \oplus \xoplus{L}^\infty
      \to
      \widehat{W}_{\widehatNeumann}^{1,q}
    }
  \end{math}
  is boundedly Lipschitzian, \ie\  for any bounded subset
  \begin{math}
    M \subset 
    \widehat{W}_{\widehatNeumann}^{-1,q}
    \oplus
    \xoplus{L}^\infty
  \end{math}
  there is a constant $\nlpsolution_M$ such that
  \begin{equation*}
    \lrnorm{
      \nlpsolution(f,z) 
      - 
      \nlpsolution(\check{f},\check{z})
    }_{W^{1,q}} 
    \le
    \nlpsolution_M
    \left(
      \lrnorm{f-\check{f} }_{\widehat{W}_{\widehatNeumann}^{-1,q}}
      +
      \lrnorm{z - \check{z} }_{\xoplus{L}^\infty}
    \right)
  \end{equation*}
  for all $(f,z)$, $(\check{f},\check{z})\in{M}$.

  2.~Let additionally \assuref{assu:dop} be satisfied. If
  \begin{equation*}
    z=(z_1,z_2) 
    \in 
    C([\TNULL,T],\xoplus{L}^\infty) \cap C^1(]\TNULL,T[,\xoplus{L}^p)
    ,
  \end{equation*}
  then the function 
  \begin{math}
    [\TNULL,T] \ni t 
    \mapsto
    \varphi(t) \in \widehat{W}_{\widehatNeumann}^{1,q}
  \end{math}
  given by
  \begin{math}
    \varphi(t) 
    \df 
    \nlpsolution(\immotile(t),z(t))
    \in 
    \widehat{W}_{\widehatNeumann}^{1,q} 
  \end{math}
  is continuous, and continuously differentiable on $]\TNULL,T[$. 
  Its derivative is 
  \begin{multline*}
    \label{eq:impl21}
    \varphi'(t) 
    =
    \left[
      \partial \nlpsolution\big(\immotile(t),z(t)\big)
    \right]
    \big(\immotile'(t),z'(t)\big)
    \\
    =\left(
      \Linpoisson + \Lindensity_1 + \Lindensity_2
    \right)^{-1}
    \left(
      \immotile'(t) 
      + \Lindensity_1  \extension z_1 ' 
      - \Lindensity_2  \extension z_2 ' 
    \right),
  \end{multline*}
  where $\Lindensity_k$ is again defined by \eqref{eq:pk} --- there
  $(F,Z)$ specified as $\big(\immotile(t),z(t)\big)$.
\end{cor}

\subsection[Derivation of the quasi-linear system]
{Derivation of the quasi-linear system}
\label{sec:der}

We start now with the reformulation of the van Roosbroeck system as
defined in \defnref{def:vanroos} as a quasi-linear parabolic system for
the continuity equations.  The aim of eliminating the electrostatic
potential in mind, we first look for a substitute for its time
derivative.  In order to achieve this, we formally differentiate
Poisson's equation \eqref{eq:poisson} with respect to time. This gives
\begin{equation}
  \label{eq:poidiff} 
  \Linpoisson \varphi'
  =  
  \immotile'
  +
  \extension (u_1' - u_2')
                                %  \quad\text{on $\widehat{\Omega}$}
  .
\end{equation}
From \eqref{eq:continuity} one obtains
\begin{equation}
  \label{eq:conteq}
  u_1' - u_2' 
  =
  \dive j_1 - \dive j_2 + 
  \reaction_1(t,\Varphi,\Fullphi) - 
  \reaction_2(t,\Varphi,\Fullphi) 
                                %  \quad \text{on $\Omega$}
  .
\end{equation}
Inserting \eqref{eq:conteq} into \eqref{eq:poidiff}, one gets 
\begin{equation}
  \label{eq:currcon}
  \Linpoisson \varphi' 
  =  \immotile' 
  + \extension 
  \big(
  \dive j_1  
  -  \dive j_2 
  +  \reaction_1(t,\Varphi,\Fullphi) 
  -  \reaction_2(t,\Varphi,\Fullphi)
  \big)  
  .
\end{equation}
Just in case, $\reaction=\reaction_1=\reaction_2$ is only
recombination, this is precisely the well known conservation law for
the total current, \cf\ \cite{gajewski93}.
Clearly, \eqref{eq:currcon} leads to
\begin{equation}
  \label{eq:currcon1}
  \restrict \varphi' 
  = 
  \restrict  
  \Linpoisson^{-1} 
  \left(
    \immotile' 
    + \extension \big(
    \dive j_1 
    - \dive j_2 
    + \reaction_1(t,\Varphi,\Fullphi) 
    - \reaction_2(t,\Varphi,\Fullphi)
    \big)
  \right) 
  .
\end{equation}
Now we differentiate \eqref{eq:car-density} (with
\eqref{eq:chemical-potential}) with respect to time and obtain
\begin{multline}
  \label{eq:diffdicht}
  u'_k 
  = 
  \rho_k 
  \mathcal{F}_k'
  (\Fullphi_k + (-1)^k\restrict \Varphi + b_k)
  \big[
  \Fullphi_k' + (-1)^k\restrict  \Varphi'  + b_k'
  \big]
  \\
  +\rho_k'
  \mathcal{F}_k
  (\Fullphi_k + (-1)^k\restrict \Varphi + b_k)
  , 
  \quad k =1,2,
\end{multline}
Pending further notice we do not write out the argument
\begin{math}
  \Fullphi_k + (-1)^k\restrict \Varphi + b_k
\end{math}
of the distribution function $\mathcal{F}_k$ and its derivative.  
We also abstain from drawing out the argument of the reaction terms
$\reaction_k$.
According to \eqref{eq:splitoff} we split
$\Varphi'=\varphi'+\dbvarphi'$ and insert \eqref{eq:diffdicht} into
the current continuity equation \eqref{eq:continuity}. Thus, we find
\begin{equation*}
  \big[
  \Fullphi_k' + (-1)^k\restrict \varphi'
  \big] 
  \rho_k 
  \mathcal{F}_k'
  - \dive j_k 
  =
  \reaction_k
  - 
  \big[
  (-1)^k \restrict \dbvarphi'+b_k'
  \big]
  \rho_k 
  \mathcal{F}_k'
  -
  \rho_k' 
  \mathcal{F}_k
  ,
  \quad  k=1,2
  .
\end{equation*}
Using \eqref{eq:currcon1} we get further
\begin{multline*}
  \rho_k 
  \mathcal{F}_k'
  \Fullphi_k'
  - \dive j_k 
  + (-1)^k 
  \rho_k 
  \mathcal{F}_k'
  \restrict 
  \Linpoisson^{-1} 
  \big(
  \immotile' 
  + 
  \extension \big(
  \dive j_1 - \dive j_2 
  + \reaction_1
  - \reaction_2
  \big)
  \big)
  \\
  = 
  \reaction_k
  - 
  \big[
  (-1)^k\restrict \dbvarphi'+ b_k'
  \big]
  \rho_k 
  \mathcal{F}_k'
  -
  \rho_k' 
  \mathcal{F}_k
  ,
  \quad  k=1,2
  .
\end{multline*}
Dividing this by $\rho_k\mathcal{F}_k'$ we obtain
\begin{multline*}
                                %  \label{eq:evol1}
  \left(
    \begin{array}{c}
      \Fullphi_1'\\
      \Fullphi_2'
    \end{array}
  \right) 
  - 
  \left(
    \begin{array}{cc}
      1 + \restrict\Linpoisson^{-1}\extension \mathcal{F}_1' \rho_1  
      & 
      -\restrict\Linpoisson^{-1}\extension \mathcal{F}_2' \rho_2
      \\
      -\restrict\Linpoisson^{-1}\extension \mathcal{F}_1' \rho_1 
      & 
      1 + \restrict\Linpoisson^{-1}\extension \mathcal{F}_2' \rho_2
    \end{array}
  \right) 
  \left(
    \begin{array}{cc}
      \frac{1}{\rho_1 \mathcal{F}_1'}
      & 
      0 
      \\
      0 
      &  
      \frac{1}{\rho_2 \mathcal{F}_2'}
    \end{array}
  \right) 
  \left(
    \begin{array}{c}
      \dive j_1\\
      \dive j_2
    \end{array}
  \right) 
  \\[1ex]                       
 = 
  \left(
    \begin{array}{c}
      \frac{\reaction_1}{\rho_1 \mathcal{F}_1'}
      +
      \reaction_1\restrict\Linpoisson^{-1}\extension 
      -
      \reaction_2\restrict\Linpoisson^{-1}\extension
      \\
      -
      \reaction_1\restrict\Linpoisson^{-1}\extension
      + 
      \frac{\reaction_2}{\rho_2 \mathcal{F}_2'}
      +
      \reaction_2\restrict\Linpoisson^{-1}\extension
    \end{array}
  \right)
%   \left(
%     \begin{array}{cc}
%       \frac{1}{\rho_1 \mathcal{F}_1'}
%       +
%       \restrict\Linpoisson^{-1}\extension 
%       & 
%       -\restrict\Linpoisson^{-1}\extension
%       \\
%       -\restrict\Linpoisson^{-1}\extension
%       & 
%       \frac{1}{\rho_2 \mathcal{F}_2'}
%       +
%       \restrict\Linpoisson^{-1}\extension
%     \end{array}
%   \right)
%   \left(
%     \begin{array}{c}
%       \reaction_1\\
%       \reaction_2
%     \end{array}
%   \right) 
  + 
  \left(
    \begin{array}{c}
      \restrict \Linpoisson^{-1} \immotile' 
      + \restrict \dbvarphi' 
      - b_1'
      -  \frac{\rho_1'}{\rho_1} \frac{\mathcal{F}_1}{\mathcal{F}'_1}
      \\
      -
      \restrict \Linpoisson^{-1} \immotile' 
      - \restrict \dbvarphi' 
      - b_2'
      - \frac{\rho_2'}{\rho_2} \frac{\mathcal{F}_2}{\mathcal{F}'_2}
    \end{array}
  \right)
\end{multline*}
This evolution equation can be written in the condensed form
\begin{equation}
  \label{eq:evol2}
  \Fullphi' - [I + Z(t,\Fullphi)]E(t,\Fullphi)\DIVE{j} = Y(t,\Fullphi)
\end{equation}
where $\Fullphi=(\Fullphi_1,\Fullphi_2)$ and
\begin{math}
  \DIVE{j} \df ( \dive j_1 , \dive j_2 ).
\end{math}
Moreover, $I$ denotes the identity.  The coefficients $Z$, $E$, and
$Y$ are given in the following way:
First we split off the Dirichlet inhomogeneities of $\Varphi$ in the
sense of \secref{sec:diribv} and we replace $\varphi$ by the solution
of the nonlinear Poisson equation, \cf\ \thmref{thm:monotone}. With
respect to an arbitrary $\psi=(\psi_1,\psi_2)\in\xoplus{W}^{1,q}$ we
set
\begin{equation}
  \label{eq:Q}
  Q_k(t,\psi) 
  \df 
  \psi_k
  + (-1)^k\restrict\nlpsolution\big(\immotile(t),z(t)\big)
  + (-1)^k\restrict\dbvarphi(t) +b_k(t)
  ,
  \quad k = 1,2,
\end{equation}
where $z\df(z_1,z_2)$ with
\begin{equation}
  \label{eq:zrhs}
  z_k (t)
  \df 
  \psi_k   
  + (-1)^k\restrict\dbvarphi(t) 
  + b_k(t)
  ,
  \quad
  k=1,2.
\end{equation}
Now we define
\begin{align}
  Z(t,\psi) 
  & \df  
  \left(
    \begin{array}{cc}
      \restrict\Linpoisson^{-1}\extension 
      \mathcal{F}_1'(Q_1(t,\psi))\rho_1(t)  
      & 
      -\restrict\Linpoisson^{-1}\extension 
      \mathcal{F}_2'(Q_2(t,\psi))\rho_2(t)
      \\
      -\restrict\Linpoisson^{-1}\extension 
      \mathcal{F}_1'(Q_1(t,\psi))\rho_1(t) 
      & 
      \restrict\Linpoisson^{-1}\extension 
      \mathcal{F}_2'(Q_2(t,\psi))\rho_2(t)
    \end{array}
  \right)
  \label{eq:Z}
  \\[1ex]
  E(t,\psi) 
  & \df 
  \left(
    \begin{smallmatrix}
      E_1(t,\psi) 
      & 
      0 
      \\
      0 
      &  
      E_2(t,\psi) 
    \end{smallmatrix}
  \right)
  ,\quad
  E_k(t,\psi) 
  \df 
  \frac{1}{\rho_k(t) \mathcal{F}_k'(Q_k(t,\psi))} 
  \label{eq:E}
  \\[1ex]
  R(t,\psi) 
  & \df 
  \left(
    \begin{array}{c}
      \reaction_1(
      t,
      \nlpsolution(\immotile(t),z(t)) + \dbvarphi(t),
      \psi)
      \\ 
      \reaction_2(
      t,
      \nlpsolution(\immotile(t),z(t)) + \dbvarphi(t),
      \psi)
    \end{array}
  \right)
  ,
  \label{eq:R}
\end{align}
and finally
\begin{equation}
  \label{eq:Y}
  Y(t,\psi) 
  \df 
  \big[
  I + Z(t,\psi)
  \big]
  E(t,\psi)R(t,\psi)
  - X(t,\psi)
  ,
\end{equation}
where $X(t,\psi)=\big(X_1(t,\psi),X_2(t,\psi)\big)$ with
\begin{equation}
  \label{eq:X}
  X_k(t,\psi) 
  \df 
  (-1)^k \restrict 
  \big(
  \Linpoisson^{-1} \immotile'(t) + \dbvarphi'(t)
  \big) 
  + b_k'(t)
  +  
  \frac{\rho_k'(t)}{\rho_k(t)} 
  \frac{\mathcal{F}_k(Q_k(t,\psi))}{\mathcal{F}'_k(Q_k(t,\psi))}
  ,\quad
  % k=1,2.
\end{equation}
$k=1,2$.
% \begin{multline}
%   \\[1ex]
%   + 
%   \left(
%     \begin{array}{c}
%       \restrict \Linpoisson^{-1} \immotile'(t) 
%       + \restrict \dbvarphi'(t) 
%       - b_1'(t)
%       -  
%       \frac{\rho_1'(t)}{\rho_1(t)} 
%       \frac{\mathcal{F}_1(Q_1(t,\psi))}{\mathcal{F}'_1(Q_1(t,\psi))}
%       \\[1ex]
%       -
%       \restrict \Linpoisson^{-1} \immotile'(t) 
%       - \restrict \dbvarphi'(t) 
%       - b_2'(t)
%       - 
%       \frac{\rho_2'(t)}{\rho_2(t)}
%       \frac{\mathcal{F}_2(Q_2(t,\psi))}{\mathcal{F}'_2(Q_2(t,\psi))}
%     \end{array}
%   \right).
% \end{multline}
Please note
\begin{equation}
  \label{eq:ZE}
  Z(t,\psi)E(t,\psi) 
  = 
  \left(
    \begin{smallmatrix}
      \restrict\Linpoisson^{-1}\extension 
      & 
      -\restrict\Linpoisson^{-1}\extension
      \\
      -\restrict\Linpoisson^{-1}\extension
      &\restrict\Linpoisson^{-1}\extension
    \end{smallmatrix}
  \right)
  .
\end{equation}

Next we apply the definition \eqref{eq:curr-dens} of the currents
$j_k$ and get
\begin{equation*}
  \dive j_k
  = 
  \dive 
  \big(
  \mathcal{G}_k(
  \Fullphi_k
  + (-1)^k \restrict \varphi 
  + (-1)^k \restrict \dbvarphi
  +b_k)
  \mu_k
  \grad \Fullphi_k 
  \big)
  , 
  \quad k = 1,2
  ,
\end{equation*}
or in shorter notation
\begin{equation}
  \label{eq:curmoda}
  \DIVE {j} 
  =
  \DIVE 
  G(t,\Fullphi)
  \mu
  \GRAD 
  \Fullphi
  ,
\end{equation}
where --- \cf\ also \eqref{eq:Q} and \eqref{eq:curr-dens} ---
\begin{equation}
  \label{eq:G}
  G(t,\psi) 
  \df
  \left(
    \begin{smallmatrix}
      G_1(t,\psi) 
      & 
      0 
      \\
      0 
      &  
      G_2(t,\psi) 
    \end{smallmatrix}
  \right)
  ,\quad
  G_k(t,\psi)
  \df 
  \mathcal{G}_k \big( Q_k(t,\psi) \big) 
  .
\end{equation}

Now, putting together \eqref{eq:curmoda} and \eqref{eq:evol2} we
obtain in conclusion the evolution equation
\begin{equation}
  \label{eq:evol3}
  \Fullphi'
                                %\frac{\partial \Fullphi}{\partial t} 
  - 
  \big[ I + Z(t,\Fullphi) \big]
  E(t,\Fullphi)
  \DIVE
  G(t,\Fullphi) \mu \GRAD \Fullphi 
  = Y(t,\Fullphi)
\end{equation}
which has to be complemented by the boundary conditions
\eqref{eq:cc:bc} and the initial condition \eqref{eq:cc:ini},
\cf\ also \remref{rem:curnormal}.

%%\section[Auxiliary evolution equation]{Auxiliary evolution equation}
\section
[The quasi-linear parabolic equation]
{The quasi-linear parabolic equation}
\label{sec:evol}

Evolution equations of the type \eqref{eq:evol3} were investigated in
\cite{pub:765}: \eqref{eq:evol3} has a unique, local in time solution,
if the functions $E$, $G$, $Z$ and $Y$ defined by \eqref{eq:E},
\eqref{eq:G}, \eqref{eq:Z} and \eqref{eq:Y}, respectively, satisfy the
following conditions.
\begin{assu}
  \label{assu:evol}
  With respect to $q\in]2,\infty[$ and $p=q/2$, as specified in
  \defnref{def:pq}, there is an $\eta\in]0,1]$ and further for any
  bounded set $M\subset\xoplus{W}^{1,q}$ exist positive constants
  $\constem$, $\constgm$, $\constym$, and $\constzm$ such that the
  mappings
  \begin{eqnarray}
    E 
                                %= \diag(E_1,E_2) 
    & : & 
    [\TNULL,\TEINS] \times \xoplus{W}^{1,q} 
    \longrightarrow 
    \xoplus{L}^\infty,
    \label{E1}\\
    G 
                                %= \diag(G_1,G_2) 
    & : & 
    [\TNULL,\TEINS] \times \xoplus{W}^{1,q} 
    \longrightarrow
    \xoplus{W}^{1,q},
    \label{G1}\\
    Z 
                                %= (Z_{km})_{k.m=1,2} 
    & : & 
    [\TNULL,\TEINS] \times \xoplus{W}^{1,q} 
    \longrightarrow 
    \mathcal{B}_\infty(\xoplus{L}^p)
    \label{Z1},\\
    Y 
                                %= (Y_1,Y_2)^T    
    & : & [\TNULL,\TEINS] \times \xoplus{W}^{1,q} 
    \longrightarrow 
    \xoplus{L}^p\label{Y1}
  \end{eqnarray}
  satisfy the conditions
  \begin{eqnarray}
    \min_{k = 1,2}
    \inf_{
      \begin{smallmatrix}
        t\in[\TNULL,\TEINS]\\
        \psi\in M
      \end{smallmatrix}
    }
    \essinf_{x\in\Omega}E_k(t,\psi)(x) & > & 0
    \label{E2}\\
    \min_{k = 1,2}
    \inf_{
      \begin{smallmatrix}
        t\in[\TNULL,\TEINS]\\
        \psi\in M
      \end{smallmatrix}
    }
    \essinf_{x\in\Omega}G_k(t,\psi)(x) & > & 0
    \label{G2}
  \end{eqnarray}
  and for all $\ta$, $\tb\in[\TNULL,\TEINS]$ and all $\psi,
  \check{\psi} \in M$:
  \begin{eqnarray}
    \|E(\ta,\psi) - E(\tb,\check{\psi})\|_{\xoplus{L}^\infty} 
    & \le &  
    \constem 
    \left(
      |\ta-\tb|^\eta 
      + \|\psi- \check{\psi}\|_{\xoplus{W}^{1,q}}
    \right), 
    \label{E3}\\
    \|G(\ta,\psi) - G(\tb,\check{\psi})\|_{\xoplus{W}^{1,q}} 
    & \le &
    \constgm 
    \left(
      |\ta-\tb|^\eta 
      + \|\psi-\check{\psi}\|_{\xoplus{W}^{1,q}}
    \right), 
    \label{G3}\\
    \|Z(\ta,\psi) - Z(\tb, \check{\psi})\|_{\mathcal{B}(\xoplus{L}^p)} 
    & \le &
    \constzm 
    \left(
      |\ta-\tb|^\eta 
      + \|\psi-\check{\psi}\|_{\xoplus{W}^{1,q}}
    \right),
    \label{Z3}\\
    \|Y(\ta,\psi) - Y(\tb, \check{\psi})\|_{\xoplus{L}^p} 
    & \le &
    \constym 
    \left(|\ta-\tb|^\eta + 
      \|\psi- \check{\psi}\|_{\xoplus{W}^{1,q}}
    \right)
    \label{Y3}
    .
  \end{eqnarray}
\end{assu}
\begin{defn}
  \label{def:solution} 
  Let the Assumptions~\ref{assu:cc:bc} and \ref{assu:evol}
  be satisfied. Further, let 
  \begin{math}
    \matrixA
    \from{
      \domainA \to \xoplus{L}^p
    }
  \end{math}
  be the operator from \defnref{def:pq} and let $V$ be a Banach space
  such that
  $\domainA\embedcontinuously{V}\embedcontinuously\xoplus{W}^{1,q}$.
  We say the evolution equation \eqref{eq:evol3} with initial
  condition $\Fullphi(\TNULL)=\iniphi\in\xoplus{W}^{1,q}$ has a unique
  local solution $\Fullphi=\phi+\bvphi{}$ with respect to $V$ if
  $\iniphi-\bvphi{}(\TNULL)\in{V}$ implies the existence of a number
  $T\in]\TNULL,\TEINS]$ such that the initial value problem
  \begin{multline}
    \label{eq:exactsys}    
    \phi'(t) 
    + 
    \big[
    I + Z\big(t,\phi(t) + \bvphi{}(t)\big)
    \big]
    E\big(t,\phi + \bvphi{}(t)\big)
    G\big(t,\phi(t) + \bvphi{}(t)\big)
    \matrixA\phi(t)
    \\
    = 
    Y\big(t,\phi(t) + \bvphi{}(t)\big) 
    - 
    \bvphi{}'(t) 
    + 
    J\big(t,\phi(t)\big),
    \quad
    \phi(\TNULL) = \iniphi - \bvphi{}(\TNULL)
  \end{multline}
  admits a unique solution
  \begin{equation}
    \label{eq:solution}
    \phi 
    \in 
    C^1(]\TNULL,T[,\xoplus{L}^p) 
    \cap 
    C(]\TNULL,T],\domainA) 
    \cap 
    C([\TNULL,T],V)
    .
  \end{equation}
  For 
  \begin{math}
    (t,\psi) \in [\TNULL,\TEINS] \times \xoplus{W}^{1,q}_{\Neumann}
  \end{math}
  the term $J$ in \eqref{eq:exactsys} is given by
  \begin{equation*}
    J(t,\psi) 
    \df 
    \big[
    I + Z\big(t,\psi + \bvphi{}(t)\big)
    \big]
    E\big(t,\psi + \bvphi{}(t)\big)
    \GRAD G\big(t,\psi + \bvphi{}(t)\big)
    \cdot 
    \mu \GRAD\big(\psi + \bvphi{}(t)\big)
    . 
  \end{equation*}
\end{defn}           
\begin{rem}
  \label{rem:equiv}
  We have to clarify the relation between \eqref{eq:evol3} and
  \eqref{eq:exactsys}.
  If
  $\Fullphi=\phi+\bvphi{}$ is a solution in the sense of
  \defnref{def:solution}, then 
  \begin{equation} 
    \label{eq:distrib}
    \DIVE 
    G(t,\Fullphi)\mu\GRAD\Fullphi
    = 
    G(t,\Fullphi)\matrixA\phi 
    +
    \GRAD G(t,\Fullphi)
    \cdot
    \mu \GRAD\Fullphi
  \end{equation}
  is satisfied, which allows to rewrite \eqref{eq:exactsys} in the form 
  \eqref{eq:evol3}.
\end{rem}
\begin{rem}
  \label{rem:boundary}
  If $\Fullphi=(\Fullphi_1,\Fullphi_2)$ is a solution of
  \eqref{eq:evol3} in the sense of \defnref{def:solution}, then
  \begin{equation*}
    \trace 
    \big(
    \Fullphi_k(t)
    \big)
    \big|_{\Dirichlet} 
    = 
    \trace 
    \big(
    \bvphi{k}(t)
    \big)
    \big|_{\Dirichlet} 
    =
    \Bvphi{k}(t) 
    ,
    \quad k = 1,2,
    \quad t \in [\TNULL,T].
  \end{equation*}
  The Neumann boundary condition 
  \begin{equation*}
    0
    =
    {\nu}\cdot{\mu_k\grad \Fullphi_k(t)}
    \big|_{\Neumann}
    = 
    {\nu}\cdot{\mu_k\grad\bvphi{k}(t)}
    \big|_{\Neumann}
    ,
    \quad k = 1,2, 
    \quad t \in [\TNULL,T],
  \end{equation*}
  holds in the distributional sense, \cf\ \remref{rem:curnormal}.
\end{rem}
\begin{prop}
  \label{prop:oint} \emph{(See \cite{pub:765}.)}
  Let the Assumptions~\ref{assu:cc:bc} and \ref{assu:evol} be
  satisfied. For each 
  \begin{math}
    \gamma 
    \in 
    \big]
    \frac{1}{2}+\frac{1}{q},1
    \big[
  \end{math}
  the initial value problem \eqref{eq:evol3} with initial value
  $\iniphi\in\xoplus{W}^{1,q}$ has a unique local
  solution $\phi$ with respect to the complex interpolation spaces
  $V\df\big[\xoplus{L}^p,\domainA\big]_\gamma$.
\end{prop}
We are now going to show that the mappings $E$, $G$, $Y$ and $Z$ satisfy 
\assuref{assu:evol}. To that end we need the following preparatory lemma.
\begin{lem} 
  \label{lem:nemyckii}
  If $\xi\from{\mathbb{R}\to\mathbb{R}}$ is continuously
  differentiable, then $\xi$ induces a Nemyckii operator from
  $L^\infty$ into itself which is boundedly Lipschitzian.
  If $\xi\from{\mathbb{R}\to\mathbb{R}}$ is twice continuously
  differentiable, then it induces a Nemyckii operator from $W^{1,q}$
  into itself which is boundedly Lipschitzian.
\end{lem}
The proof is straightforward. Recall that, according to
\defnref{def:pq}, $q$ is fixed and larger than two.
\begin{lem} 
  \label{lem:contQ}
  Let the Assumptions~\ref{assu:bands}, \ref{assu:poi-diri-bv} and
  \ref{assu:dop} be satisfied. Then the equation \eqref{eq:Q} defines
  mappings 
  \begin{math}
    Q_k \from{ [\TNULL,\TEINS] \times \xoplus{L}^\infty 
      \to
      L^\infty
    }, 
  \end{math}
  $k =1,2$, and the restriction of each $Q_k$ to
  $[\TNULL,\TEINS]\times\xoplus{W}^{1,q}$ takes its values in
  $W^{1,q}$. Moreover, there is a number $\eta\in]0,1]$ and then for
  any bounded subset $M\subset\xoplus{L}^\infty$ a positive number
  $\constqm$ exists such that
  for all $\ta,\,\tb\in[\TNULL,\TEINS]$ and all
  $\psi,\,\check{\psi}\in M$:
  \begin{equation*}
    \|
    Q_k(\ta,\psi) - Q_k(\tb,\check{\psi})
    \|_{L^\infty} 
    \le
    \constqm
    \big(
    |\ta - \tb|^\eta 
    +\|\psi - \check{\psi}\|_{\xoplus{L}^\infty}
    \big),
    \quad k = 1,2.
  \end{equation*}
  Analogously, for each bounded subset $M\subset\xoplus{W}^{1,q}$
  there is a positive number $\constqm$ such that for all
  $\ta,\,\tb\in[\TNULL,\TEINS]$ and all $\psi,\,\check{\psi}\in M$:
  \begin{equation*}
    \|
    Q_k(\ta,\psi) - Q_k(\tb,\check{\psi})
    \|_{W^{1,q}} 
    \le 
    \constqm 
    \big(
    |\ta-\tb|^\eta 
    + \|\psi - \check{\psi}\|_{\xoplus{W}^{1,q}}
    \big),
    \quad k = 1,2.
  \end{equation*}
\end{lem}
The proof is obtained from \corref{cor:boundedlip}.
\begin{lem} 
  \label{lem:comp} 
  Let the Assumptions~\ref{assu:bands}, \ref{assu:poi-diri-bv} and
  \ref{assu:dop} be satisfied.
  If $\xi\from{\mathbb{R}\to\mathbb{R}}$ is continuously
  differentiable, then $\xi$ induces operators
  \begin{equation*}
    [\TNULL,\TEINS] \times \xoplus{L}^\infty 
    \ni 
    (t,\psi) 
    \longmapsto 
    \xi(Q_k(t,\psi))
    \in 
    L^\infty,
    \quad
    k = 1,2.
  \end{equation*}
  Moreover, there is a constant $\eta\in]0,1]$ and for any bounded set
  $M\subset\xoplus{L}^\infty$ a constant $\xi_M$ such that for all
  $\ta,\,\tb\in[\TNULL,\TEINS]$ and all $\psi,\,\check{\psi}\in M$:                                %
  \begin{equation*}
    \|
    \xi\big(Q_k(\ta,\psi)\big) 
    - 
    \xi\big(Q_k(\tb,\check{\psi})
    \big) 
    \|_{L^\infty} 
    \le 
    \xi_M
    \big(
    |\ta - \tb|^\eta 
    + 
    \|\psi - \check{\psi}\|_{\xoplus{L}^\infty}
    \big),
    \quad
    k=1,2.
  \end{equation*}
  If $\xi$ is twice continuously differentiable, then the restriction
  of $\xi{\circ}Q_k$ to $[\TNULL,\TEINS]\times\xoplus{W}^{1,q}$ maps
  into $W^{1,q}$, $k =1,2$. Moreover, there is a number $\eta\in]0,1]$
  and for any bounded subset $M\subset\xoplus{W}^{1,q}$ a constant
  $\xi_M$ such that for all $\ta,\,\tb\in[\TNULL,\TEINS]$ and all
  $\psi,\,\check{\psi}\in M$:
  \begin{equation*}
    \|
    \xi\big(Q_k(\ta,\psi)\big) 
    - 
    \xi\big(Q_k(\tb,\check{\psi})\big) 
    \|_{W^{1,q}} 
    \le 
    \xi_M
    \big(
    |\ta - \tb|^\eta 
    + 
    \|\psi - \check{\psi}\|_{\xoplus{W}^{1,q}}
    \big), 
    \quad
    k=1,2.
  \end{equation*}
\end{lem}
The proof follows from \lemref{lem:nemyckii} and \lemref{lem:contQ}.
\begin{lem}
  \label{lem:EG}
  Let the Assumptions~\ref{assu:bands}, \ref{assu:poi-diri-bv} and
  \ref{assu:dop} be satisfied.  Then there is a number $\eta\in]0,1]$
  such that the mappings $E$ and $G$ defined by \eqref{eq:E} and
  \eqref{eq:G} satisfy the conditions \eqref{E1}, \eqref{E2},
  \eqref{E3}, and \eqref{G1}, \eqref{G2}, \eqref{G3}, respectively.
\end{lem}
\begin{proof}
  The functions $\frac{1}{\mathcal{F}'_k}$ are continuously
  differentiable by \assuref{assu:distri}. Consequently, by
  \lemref{lem:comp} the mappings $\widetilde{E}_k$, given by
  \begin{equation*}
    [\TNULL,\TEINS] \times \xoplus{L}^\infty 
    \ni 
    (t,\psi) 
    \longmapsto
    \frac{1}{\mathcal{F}'_k\big(Q_k(t,\psi)\big)}
    \in
    L^\infty
    ,
    \quad
    k =1,2,
  \end{equation*}
  are well defined.
  Moreover, \lemref{lem:comp} provides a constant $\eta\in]0,1]$ such
  that for any bounded set $M\subset\xoplus{L}^\infty$ a constant
  $C_M$ exists such that for all
  $\ta,\,\tb\in[\TNULL,\TEINS]$ and all $\psi,\,\check{\psi}\in M$:
  \begin{equation*}
    \|
    \widetilde{E}_k(\ta,\psi) 
    - 
    \widetilde{E}_k(\tb,\check{\psi})
    \|_{\xoplus{L}^\infty} \le
    C_M
    \big(
    |\ta - \tb|^\eta 
    + 
    \|\psi - \check{\psi}\|_{\xoplus{L}^\infty}
    \big),
    \quad
    k=1,2.
  \end{equation*}
  Since $\xoplus{W}^{1,q}$ embeds continuously into
  $\xoplus{L}^\infty$ for any bounded set $M\subset\xoplus{W}^{1,q}$
  there is a constant, again named $C_M$, such that for all
  $\ta,\,\tb\in[\TNULL,\TEINS]$ and all $\psi,\,\check{\psi}\in M$:
  \begin{equation*}
    \|
    \widetilde{E}_k(\ta,\psi) 
    - 
    \widetilde{E}_k(\tb,\check{\psi})
    \|_{\xoplus{L}^\infty} 
    \le
    C_M
    \big(
    |\ta - \tb|^\eta 
    + 
    \|\psi - \check{\psi}\|_{\xoplus{W}^{1,q}}
    \big),
    \quad
    k=1,2.
  \end{equation*}
  The identity $E_k=\frac{1}{\rho_k}\widetilde{E}_k$ and
  \assuref{assu:effband} now imply \eqref{E1} and \eqref{E3}. According to 
  \lemref{lem:contQ} the sets
  \begin{equation*}
    \left\{
      Q_k(t,\phi)
      \with
      (t,\phi) \in [\TNULL,\TEINS] \times M
    \right\}
    ,
    \quad 
    k =1,2,
  \end{equation*}
  are bounded in $L^\infty$. Since the derivative of the carrier
  distribution functions $\mathcal{F}_k$, $k =1,2$, are continuous and
  positive, \eqref{E2} immediately follows.
  
  Using the second assertion of \lemref{lem:comp} we verify
  \eqref{G1}, \eqref{G2}, and \eqref{G3} in a similar manner.
\end{proof}
\begin{lem}
  \label{lem:BZ}
  Let the Assumptions~\ref{assu:bands}, \ref{assu:poi-diri-bv}, and
  \ref{assu:dop} be satisfied.  Then the mapping $Z$ given by
  \eqref{eq:Z} defines a family 
  \begin{math}
    \{Z(t,\psi)\}_{(t,\psi) \in [\TNULL,\TEINS] \times \xoplus{W}^{1,q}}
  \end{math}
  of linear, compact operators
  \begin{math}
    Z(t,\phi)\from{\xoplus{L}^p \to \xoplus{L}^p}
    .
  \end{math}
  Additionally, there is a H\"older exponent $\eta\in]0,1]$ and
  constants $\constzm$ such that \eqref{Z1} and \eqref{Z3} are
  satisfied.
\end{lem}
\begin{proof}
  It suffices to show the analogous assertions for the entries of the
  operator matrices $Z(t,\psi)$.  Firstly, \lemref{lem:comp} gives us
  the estimate
  \begin{multline*} 
    \| 
    \mathcal{F}_k'\big(Q_k(\ta,\psi)\big) 
    -
    \mathcal{F}_k'\big(Q_k(\tb,\check{\psi})\big)
    \|_{\mathcal{B}(L^p)} 
    \\
    \le
    \| 
    \mathcal{F}_k'\big(Q_k(\ta,\psi)\big) 
    -
    \mathcal{F}_k'\big(Q_k(\tb,\check{\psi})\big)
    \|_{L^\infty} 
    \\
    \le 
    C_M 
    \big(
    |\ta-\tb|^\eta 
    +
    \|\psi-\check{\psi}\|_{\xoplus{W}^{1,q}} 
    \big)
    ,
    \quad 
    k =1,2,
  \end{multline*}
  where the constant $C_M$ can be taken uniformly with respect to
  $\ta,\,\tb\in[\TNULL,\TEINS]$ and $\psi,\,\check{\psi}$ from any
  bounded set $M\subset\xoplus{W}^{1,q}$. This estimate together with
  \assuref{assu:effband} implies \eqref{Z3}. As
  $\restrict\Linpoisson^{-1}\extension$ is a linear and even compact
  operator from $L^p$ into itself, this gives \eqref{Z1}.
\end{proof}
\begin{lem}
  \label{lem:YR}
  Let the Assumptions~\ref{assu:recomb}, \ref{assu:bands},
  \ref{assu:poi-diri-bv}, and \ref{assu:dop} be satisfied.  Then the
  mapping $Y$ defined by \eqref{eq:Y} meets the conditions
  \eqref{Y1} and \eqref{Y3}.
\end{lem}
\begin{proof}
  At first one deduces from the assumptions and
  \corref{cor:boundedlip} that \eqref{eq:R} defines a mapping
  \begin{math}
    R \from{
      [\TNULL,\TEINS] \times \xoplus{W}^{1,q} \to \xoplus{L}^p
    }
  \end{math}
  for which there is a H\"older exponent $\eta\in]0,1]$. Moreover, for
  any bounded set $M\subset\xoplus{W}^{1,q}$ exists a constant $C_M$
  such that for all $\ta,\,\tb\in[\TNULL,\TEINS]$ and all
  $\psi,\,\check{\psi}\in M$:
  \begin{equation*}
    \|
    R(\ta,{\psi}) 
    - 
    R(\tb,\check{\psi})
    \|_{\xoplus{L}^p} 
    \le 
    C_M
    \big(
    |\ta - \tb|^\eta
    + 
    \|{\psi} - \check{\psi}\|_{\xoplus{W}^{1,q}}
    \big)
    .    
  \end{equation*}
  Applying \lemref{lem:EG} and \lemref{lem:BZ} one obtains \eqref{Y1}
  and \eqref{Y3} for the mapping
  \begin{equation*}
    [\TNULL,\TEINS] \times  \xoplus{W}^{1,q}
    \ni
    (t,\psi) 
    \longmapsto 
    \big[
    I + Z(t,\psi)
    \big]
    E(t,\psi)
    R(t,\psi).
  \end{equation*}
  The addends $b_k'$ and $\restrict\dbvarphi'$ of \eqref{eq:X} have the
  required properties due to \assuref{assu:bands} and
  \assuref{assu:poi-diri-bv}, respectively.
  For $\Linpoisson^{-1}\immotile'$ they follow from
  \assuref{assu:poi-diri-bv} (\cf\ also \remref{r-extend}),
  \assuref{assu:dop} and the fact that $\Linpoisson$ is an isomorphism
  from $\widehat{W}^{1,q}_{\widehatNeumann}$ onto
  $\widehat{W}^{-1,q}_{\widehatNeumann}$.
  The addend 
  \begin{math}
    \frac{\rho_k'(t)}{\rho_k(t)} 
    \frac{\mathcal{F}_k(Q_k(t,\psi))}{\mathcal{F}'_k(Q_k(t,\psi))}
  \end{math}
  of \eqref{eq:X} can be treated by means of \lemref{lem:comp} and
  \assuref{assu:effband}.
\end{proof}
We are now going to establish existence and uniqueness of a local
solution to the evolution equation \eqref{eq:evol3}.
\begin{thm}
  \label{thm:existenceA}
  Under the Assumptions~\ref{assu:recomb}, \ref{assu:bands},
  \ref{assu:cc:bc}, \ref{assu:cc:ini}, \ref{assu:poi-diri-bv} and
  \ref{assu:dop} the quasi-linear parabolic equation \eqref{eq:evol3}
  with the initial condition $\Fullphi(\TNULL)=\iniphi$ admits a
  unique local solution in the sense of \defnref{def:solution} with
  respect to the interpolation space
  $V=[\xoplus{L}^p,\domainA]_\theta$.
\end{thm}
\begin{proof}
  According to the Lemmas~\ref{lem:EG}, \ref{lem:BZ}, \ref{lem:YR} the
  mappings $E$, $G$, $Z$, and $Y$, defined by \eqref{eq:E},
  \eqref{eq:G}, \eqref{eq:Z}, and \eqref{eq:Y}, respectively, fulfill
  \assuref{assu:evol}.
                                %the conditions \eqref{E1}--\eqref{Y3}. 
  Hence, the result follows from \proref{prop:oint}, see also
  Remarks~\ref{rem:equiv} and \ref{rem:boundary}.
\end{proof}

% \section[Justification of the reformulation]
% {Justification of the reformulation}
% \label{s-justification}

\section[Main result]
{Main result}
\label{s-main}

We are going to show that a solution of the evolution equation
\eqref{eq:evol3} in the sense of \defnref{def:solution}
provides a solution of the van Roosbroeck system in the sense of
\defnref{def:vanroos}.

We start with a technical lemma.
\begin{lem}
  \label{lem:diff}
  Let $\xi\from{\mathbb{R}\to\mathbb{R}}$ be twice
  continuously differentiable. 
  The composition $\xi\circ\psi$ is from
  $C([\TNULL,T],L^\infty)$, if $\psi{\in}C([\TNULL,T],L^\infty)$.
  If $\psi$ composed with the embedding
  $L^\infty{\embedcontinuously}L^p$, $p\ge1$, is continuously
  differentiable in $L^p$ on $]\TNULL,T[$, then $\xi\circ\psi$
  composed with the same embedding is continuously differentiable in
  $L^p$ on $]\TNULL,T[$ and its derivative is given by
  \begin{equation*} 
    %\label{eq:diffbar}
    \frac{d\xi\circ\psi}{dt}(t)
    = \xi'\big(\psi(t)\big)\psi'(t) \in L^p 
    ,
    \quad
    t \in ]\TNULL,T[.
  \end{equation*}
\end{lem}
\begin{proof}
  If $h_1$, $h_2\in L^\infty$, then, by \lemref{lem:cardens} ---
  see also \assuref{assu:distri}, we may write
  \begin{equation*}
    \xi(h_1) - \xi(h_2)
    =
    \xi'(h_1)(h_1 - h_2)
    +
    T(h_1,h_2)((h_1 - h_2)
  \end{equation*}
  where $T(h_1,h_2)$ converges to zero in $L^\infty$ if
  $h_1{\in}L^\infty$ is fixed and $h_2$ approaches $h_1$ in the
  $L^\infty$-norm. Now we set $h_1=\psi(\ta)$ and $h_2=\psi(\tb)$ and
  divide both sides by $\ta-\tb$. In the limit $\tb\to\ta$ there is
  $\lim_{\tb\to\ta}T(\psi(\ta),\psi(\tb))=0$ in $L^\infty$, while
  $\lim_{\tb\to\ta}\frac{\psi(\ta)-\psi(\tb)}{\ta-\tb}=\psi'(\ta)$ in
  $L^p$ by supposition.
\end{proof}
Our next aim is to justify formula \eqref{eq:diffdicht}.
\begin{lem}
  \label{lem:implic}
  Let the Assumptions~\ref{assu:bands}, \ref{assu:cc:bc},
  \ref{assu:poi-diri-bv}, and \ref{assu:dop} be satisfied and assume that
  $\Fullphi$ is a solution of \eqref{eq:evol3}.  We define
  \begin{equation}
    \label{eq:zimpk}
    z\df(z_1,z_2)
    \quad\text{with}\quad
    z_k(t) \df \Fullphi_k(t) + b_k(t) + (-1)^k \restrict \dbvarphi(t)
    ,
    \quad
    k=1,2
    ,
    \;
    t \in [\TNULL,T]
    ,
  \end{equation}
  and $\varphi(t)\df\nlpsolution\big(\immotile(t),z(t)\big)$.  Then
  \begin{math}
    Q_k(t,\Fullphi(t)) 
    =
    z_k(t) 
    +
    (-1)^k \restrict\varphi(t)
    ,
  \end{math}
  and
  the functions
  \begin{equation*}
    [\TNULL,T] \ni t 
    \longmapsto 
    G_k(t,\Fullphi(t))
    =
    \mathcal{G}_k 
    \big(
    Q_k(t,\Fullphi(t)) 
    \big) \in L^\infty 
    ,
  \end{equation*}
  and
  \begin{equation*} 
    % \label{eq:uk}
    [\TNULL,T] \ni t 
    \longmapsto 
    u_k(t) \df \rho_k(t) \mathcal{F}_k
    \big( Q_k(t,\Fullphi(t)) 
    \big) \in L^\infty 
  \end{equation*}
  are continuous and concatenated with the embedding
  $L^\infty{\embedcontinuously}L^p$ they are continuously
  differentiable on $]\TNULL,T[$. The time derivative of
  $u_k$ %\eqref{eq:uk}
  is given by
  \begin{multline} 
    \label{eq:ukprime}
    u_k'(t)  
    = 
    \rho_k'(t) \mathcal{F}_k
    \big( 
    Q_k(t,\Fullphi(t))
    % z_k(t) + (-1)^k \restrict \varphi(t) 
    \big)
    \\  
    +
    \rho_k(t) \mathcal{F}_k'
    \big( 
    Q_k(t,\Fullphi(t))
    % z_k(t) + (-1)^k \restrict \varphi(t) 
    \big)
    \big[
    \Fullphi_k'(t) 
    + b_k'(t) 
    + (-1)^k \restrict \dbvarphi'(t)
    + (-1)^k \restrict \varphi'(t)
    \big]
  \end{multline}
  $k =1,2$, $t\in]\TNULL,T]$.
\end{lem}
\begin{proof}
  Due to \assuref{assu:cc:bc} and \defnref{def:solution} the function
  $\Fullphi$ belongs to the space 
  \begin{equation} 
    \label{eq:space}
    C([\TNULL,T],\xoplus{L}^\infty) 
    \cap 
    C^1(]\TNULL,T[,\xoplus{L}^p)
  \end{equation} 
  see also \remref{rem:initial}.
  Hence, the Assumptions~\ref{assu:bands} and \ref{assu:poi-diri-bv}
  ensure that the function $z$ also belongs to this space, and by
  \corref{cor:boundedlip}, so does the function
  $\varphi=\nlpsolution\big(\immotile(t),z(t)\big)$.  Thus, we may
  apply \lemref{lem:diff}.
\end{proof}
\begin{rem} 
  \label{rem:justi}
  \lemref{lem:implic} justifies the formal manipulations in
  \secref{sec:der}. First, \eqref{eq:diffdicht} is given a strict
  sense. Furthermore, the differentiation of Poisson's equation
  \eqref{eq:poidiff} has the following precise meaning: since
  $\Fullphi$ is from the space \eqref{eq:space}, the function
  $t\mapsto\varphi(t)$ is differentiable --- even in a much 'better'
  space than $\Fullphi$ --- \cf\ \corref{cor:boundedlip}. Hence, the
  right hand side of \eqref{eq:poisson} is differentiable with respect
  to time in the space $\widehat{W}^{-1,q}_{\widehatNeumann}$ and
  \eqref{eq:poidiff} is an equation in the space
  $\widehat{W}^{-1,q}_{\widehatNeumann}$.
\end{rem}
We come now to the main results of this paper.
\begin{thm} 
  \label{thm:central}
  Under the Assumptions~\ref{assu:recomb}, \ref{assu:bands},
  \ref{assu:cc:bc}, \ref{assu:cc:ini}, \ref{assu:poi-diri-bv}, and
  \ref{assu:dop} van Roosbroeck's system with initial condition
  $\Fullphi(\TNULL)=\iniphi\in\xoplus{W}^{1,q}$ admits a unique local
  in time solution in the sense of \defnref{def:vanroos}.
\end{thm}
\begin{proof} 
  By \thmref{thm:existenceA} the auxiliary evolution equation
  \eqref{eq:evol3} admits --- in the sense of \defnref{def:solution}
  --- a unique local solution $\Fullphi$ satisfying the initial
  condition $\Fullphi(\TNULL)=\iniphi$.  Let us show that --- in the sense
  of \defnref{def:vanroos} --- the pair $\{\Varphi,\Fullphi\}$, with
  $\Varphi$ given by
  \begin{equation}
    \label{eq:proof1}
    \Varphi(t) 
    \df  
    \dbvarphi(t)
    +
    \nlpsolution
    \big(
    \immotile(t),
    z(t)
    \big) 
    , 
    \quad t \in [\TNULL,T],
  \end{equation}
  and $z$ according to \eqref{eq:zimpk},
  is a local solution of van Roosbroeck's system. 
  First, \eqref{eq:sol-phi} is identical with \eqref{eq:solution}.  By
  the embedding 
  \begin{math}
    V 
    \embedcontinuously 
    \xoplus W^{1,q}_{\Neumann}
    \embedcontinuously 
    \xoplus{L}^\infty
  \end{math}
  (\cf\ \remref{rem:initial}) the
  function 
  \begin{math}
    [\TNULL,T]\ni{t}\mapsto\phi(t)\in \xoplus{L}^\infty
  \end{math}
  is continuous, and so is the function 
  \begin{math}
    [\TNULL,T]\ni{t}\mapsto\bvphi{}(t)\in  \xoplus{L}^\infty
  \end{math}
  in view of \assuref{assu:cc:bc}.
  Thus, 
  \begin{math}
    \Fullphi 
    \in 
    C([\TNULL,T],\xoplus{L}^\infty) 
    \cap
    C^1(]\TNULL,T[,\xoplus{L}^p)
    .
  \end{math}
  Moreover, for $z$, \cf\ \eqref{eq:zimpk}, one obtains from the
  Assumptions~\ref{assu:bands} and \ref{assu:poi-diri-bv} that
  \begin{math}
    z 
    \in C([\TNULL,T],\xoplus{L}^\infty) 
    \cap
    C^1(]\TNULL,T[,\xoplus{L}^p)
    .
  \end{math}
  Consequently, property \eqref{eq:sol-varphi} follows by
  \corref{cor:boundedlip}, while \eqref{eq:ukjk-regularity} results from
  \lemref{lem:implic}. The Poisson equation \eqref{eq:poisson} with
  densities \eqref{eq:car-density} is obviously satisfied by
  \eqref{eq:proof1} due to the definition of $\nlpsolution$.
  \eqref{eq:jk-regularity} follows from 
  \begin{math}
    \grad \Fullphi_k 
    \in
    C(]\TNULL,T],\xoplus{L}^q)
    ,
  \end{math}
  $k =1,2$, and \lemref{lem:implic}.
  \eqref{eq:divjk-regularity} is implied by \eqref{eq:solution} and
  \eqref{eq:distrib}.  
  It remains to show that the continuity equations
  \eqref{eq:continuity} are satisfied.  
  For this, one first notes the relations
  \begin{equation}
    \label{eq:qk}
    Q_k(t,\Fullphi(t)) 
    = \Fullphi_k(t)
    + (-1)^k\restrict\Varphi(t)
    + b_k(t)
    = z_k(t) 
    + (-1)^k\restrict\varphi(t)
    ,
    \quad k = 1,2,
  \end{equation}
  and
  \begin{equation} 
    \label{eq:recover}
    R(t,\Fullphi(t))
    =
    \left( 
      \begin{smallmatrix}
        \reaction_1(t, \Varphi(t), \Fullphi(t)) 
        \\
        \reaction_2(t, \Varphi(t), \Fullphi(t))
      \end{smallmatrix}
    \right)
    ,
  \end{equation}
  which follows from the definitions \eqref{eq:Q} and \eqref{eq:R} of
  $R$ and $Q$, and \eqref{eq:zimpk}, \eqref{eq:proof1}.
  Further, in \assuref{assu:recomb} we demand that the mappings
  $\reaction_k$, $k=1,2$, take their values in $L^p$ --- consequently,
  $R$ takes its values in $\xoplus{L}^p$.
  From \eqref{eq:ukprime} and \eqref{eq:E} one gets
  \begin{equation*}
    E_k(t,\Fullphi(t)) u_k'(t)
    = \Fullphi_k'(t) 
    + b_k'(t) 
    + (-1)^k \restrict \Varphi'(t)  
    + \tfrac{\rho_k'(t)}{\rho_k(t)}
    \tfrac{\mathcal{F}_k ( Q_k(t,\Fullphi(t) ) )}
    {\mathcal{F}_k' ( Q_k(t,\Fullphi(t) ) )}
    ,
  \end{equation*}
  and by means of the evolution equation \eqref{eq:evol3} we obtain
  \begin{multline*}
    E(t,\Fullphi(t)) u'(t)
    =
    \big[ I + Z(t,\Fullphi(t)) \big]
    E(t,\Fullphi(t))
    \DIVE
    G(t,\Fullphi(t)) \mu \GRAD \Fullphi(t) 
    \\
    +
    \big[
    I + Z(t,\Fullphi(t))
    \big]
    E(t,\Fullphi(t))
    R(t,\Fullphi(t))
    + 
    \left(
      \begin{smallmatrix}
        \restrict\Linpoisson^{-1} \immotile'(t) - \restrict\varphi'(t)
        \\
        \restrict\varphi'(t) - \restrict\Linpoisson^{-1} \immotile'(t)  
      \end{smallmatrix}
    \right)
    .
  \end{multline*}
  We now make use of the representation \eqref{eq:cur-density} of the
  currents $j=(j_1,j_2)$, and get
  \begin{multline*}
    E(t,\Fullphi(t)) 
    \left[
      u'(t) - \DIVE j(t) - R(t,\Fullphi(t))
    \right]
    \\
    =
    Z(t,\Fullphi(t))
    E(t,\Fullphi(t))
    \left[
      \DIVE j(t) + R(t,\Fullphi(t))
    \right]
    +
%     \left(
%       \begin{smallmatrix}
%         -1 & 0 \\ 0 & 1
%       \end{smallmatrix}
%     \right)
    \left(
      \begin{smallmatrix}
        \restrict\Linpoisson^{-1} \immotile'(t) - \restrict\varphi'(t)
        \\
        \restrict\varphi'(t) - \restrict\Linpoisson^{-1} \immotile'(t)  
      \end{smallmatrix}
    \right)
    .
  \end{multline*}
  We already know that the formal differentiation of Poisson's
  equation is justified, \cf\ \remref{rem:justi}. Thus,
  \eqref{eq:poidiff} yields
  \begin{multline*}
    E(t,\Fullphi(t)) 
    \left[
      u'(t) - \DIVE j(t) - R(t,\Fullphi(t))
    \right]
    \\
    =
    Z(t,\Fullphi(t))
    E(t,\Fullphi(t))
    \left[
      \DIVE j(t) + R(t,\Fullphi(t))
    \right]
    +
    \left(
      \begin{smallmatrix}
        \restrict\Linpoisson^{-1}\extension
        (u_2'(t) - u_1'(t))  
        \\
        \restrict\Linpoisson^{-1}\extension
        (u_1'(t) - u_2'(t))  
      \end{smallmatrix}
    \right)
    ,
  \end{multline*}
  and, observing \eqref{eq:ZE} and \eqref{eq:recover}, we get
  \begin{equation}
    \label{eq:proof2}
    \left[
      E(t,\Fullphi(t))
      +
      \left(
        \begin{smallmatrix}
          \restrict\Linpoisson^{-1}\extension 
          & 
          -\restrict\Linpoisson^{-1}\extension
          \\
          -\restrict\Linpoisson^{-1}\extension
          &\restrict\Linpoisson^{-1}\extension
        \end{smallmatrix}
      \right)
    \right]
    \left(
      \begin{smallmatrix}
        u_1'(t) - \dive j_1(t) - \reaction_1(t, \Varphi(t), \Fullphi(t)) 
        \\
        u_2'(t) - \dive j_2(t) - \reaction_2(t, \Varphi(t), \Fullphi(t)) 
      \end{smallmatrix}
    \right)
    =
    0
    .
  \end{equation}
  The operator on the left is continuous on $\xoplus{L}^p$; we show now
  that its kernel is trivial. Let $f_1$, $f_2\in{L}^p$ be such that
  \begin{equation*}
    \left[
      E(t,\Fullphi(t))
      +
      \left(
        \begin{smallmatrix}
          \restrict\Linpoisson^{-1}\extension 
          & 
          -\restrict\Linpoisson^{-1}\extension
          \\
          -\restrict\Linpoisson^{-1}\extension
          &\restrict\Linpoisson^{-1}\extension
        \end{smallmatrix}
      \right)
    \right]
    \left(
      \begin{smallmatrix}
        f_1 \\ f_2
      \end{smallmatrix}
    \right)
    = 0
    .
    % \quad \text{for some $f_1,f_2 \in L^p$}
  \end{equation*}
  This is equivalent to the relations
  \begin{equation*}
    f_2=-\tfrac{E_1(t,\Fullphi(t))} {E_2(t,\Fullphi(t))}f_1
    \quad\text{and}\quad
    \restrict \Linpoisson^{-1} \extension 
    \left(
      \big(
      1+\tfrac {E_1(t,\Fullphi(t))}{E_2(t,\Fullphi(t))}
      \big)
      f_1
    \right)
    =
    -E_1(t,\Fullphi(t))f_1
    .
  \end{equation*}
  \begin{math}
    \Linpoisson^{-1}\extension
    \big(
    (1+\frac {E_1(t,\Fullphi(t))}{E_2(t,\Fullphi(t))})
    f_1
    \big)
  \end{math}
  is a continuous mapping from ${W}^{1,q}_{\Neumann}$
  into $\widehat{L}^\infty$. Indeed, the embedding
  $\widehat{L}^p\embedcontinuously\widehat{W}^{-1,q}_{\widehatNeumann}$
  is continuous, and $\Linpoisson$ is an isomorphism between
  $\widehat{W}^{1,q}_{\widehatNeumann}$ and
  $\widehat{W}^{-1,q}_{\widehatNeumann}$, \cf\
  \proref{prop:isomorphy}.
  Hence, we may multiply both sides with
  \begin{math}
    f_1+\frac{E_1(t,\Fullphi(t))}{E_2(t,\Fullphi(t))}f_1
  \end{math}
  and integrate over 
  $\Omega$; this yields
  \begin{multline}
    \label{eq:glei}
    \int_{\Omega} 
    \restrict \Linpoisson^{-1} \extension 
    \left(
      f_1 
      + 
      \tfrac{E_1(t,\Fullphi(t))}{E_2(t,\Fullphi(t))} f_1 
    \right) 
    \left(
      f_1+\tfrac{E_1(t,\Fullphi(t))}{E_2(t,\Fullphi(t))} f_1
    \right)
    \der{x}
    \\
    =
    \int_{\widehat{\Omega}} 
    \Linpoisson^{-1} 
    \extension 
    \left(
      f_1 + \tfrac{E_1(t,\Fullphi(t))}{E_2(t,\Fullphi(t))} f_1
    \right) 
    \extension 
    \left(
      f_1 + \tfrac{E_1(t,\Fullphi(t))}{E_2(t,\Fullphi(t))} f_1
    \right) 
    \der{x}
    \\
    =
    -
    \int_{\Omega} 
    E_1(t,\Fullphi(t))
    \left(
      1+\tfrac{E_1(t,\Fullphi(t))}{E_2(t,\Fullphi(t))}
    \right)
    f_1^2
    \der{x}
  \end{multline}
  The quadratic form 
  \begin{math}
    \psi
    \mapsto
    \int_{\widehat{\Omega}}
    (\Linpoisson^{-1} \psi) \psi \der{x}
  \end{math}
  is non-negative on $\widehat{L}^2$ and extends by continuity to
  $\widehat{L}^p$, where it is also non-negative. On the other hand,
  the function
  \begin{math}
    E_1(t,\Fullphi(t))
    \left(
      1+\tfrac{E_1(t,\Fullphi(t))}{E_2(t,\Fullphi(t))}
    \right) 
  \end{math}
  is almost everywhere on $\Omega$ strictly positive. Therefore, the
  right hand side of \eqref{eq:glei} can only be non-negative if $f_1$
  is zero almost everywhere on $\Omega$. 
  Hence, \eqref{eq:proof2} establishes the continuity equations
  \eqref{eq:continuity}.

  To prove uniqueness of a solution of van Roosbroeck's system in the
  sense of \defnref{def:vanroos} one assures that any solution in the
  sense of \defnref{def:vanroos} procures a solution in the sense of
  \defnref{def:solution}. Indeed this has been done on a formal stage
  by the reformulation of van Roosbroeck's system as a quasi-linear
  parabolic system in \secref{sec:nl-reform}.  In fact, all formal
  steps can be carried out in the underlying function spaces.  We
  accomplish this in the sequel for the crucial points.
  \eqref{eq:poisson} and \eqref{eq:car-density} ensure, that $\varphi$
  is a solution of \eqref{eq:nlp}.  Hence, \corref{cor:boundedlip}
  implies that $\varphi$ indeed is continuously differentiable in
  $\widehat{W}_{\widehatNeumann}^{1,q}$, and, consequently,
  \eqref{eq:currcon} makes sense in
  $\widehat{W}_{\widehatNeumann}^{-1,q}$.  The derivation of
  \eqref{eq:car-density}, \cf\ also \eqref{eq:chemical-potential}, is
  justified by \lemref{lem:diff}.  Thus, \eqref{eq:diffdicht} holds in
  a strict sense.  The division by $\rho_k\mathcal{F}_k'$ is allowed
  because both factors have (uniform) upper and lower bounds. The rest
  of the manipulations up to \eqref{eq:evol3} is straight forward to
  justify.
\end{proof}
Next we want to establish the natural formulation of the balance laws
in van Roosbroeck's system in integral form, \cf\
\eqref{eq:balancelaw}, which is one of the central goals of this
paper.
At first, one realizes that the boundary integral has to be understood
in the distributional sense --- as is well known from Navier-Stokes
theory, see \cite{temam} --- if one only knows that the current is a
$q$--summable function and that its divergence is $p$--summable. More precisely, the following proposition holds.
\begin{prop} 
  \label{prop:currnor}
  Let $\subomega\subset\mathbb{R}^2$ be any bounded Lipschitz domain.
  Assume $j\from{\subomega\to\mathbb{R}^2}$ to be from
  $L^q(\subomega;\mathbb{R}^2)$ and let the divergence (in the sense
  of distributions) $\dive{j}$ of $j$ be $p$--integrable on $\subomega$.
  If $q>2$ and $p=\frac{q}{2}$, then there
  is a uniquely determined linear continuous functional 
  \begin{math}
    j_\nu 
    \in
    W^{-1+\frac{1}{q'},q}(\partial\subomega)
  \end{math}
  such that
  \begin{equation}
    \label{eq:currnor}
    \int_\subomega  j \cdot \grad\psi \der{x}
    + 
    \int_\subomega \psi \dive j \der{x}
    = 
    \dual{j_\nu}{\psi|_{\partial\subomega}}
    \quad 
    \text{for all $\psi \in W^{1,q'}(\subomega)$,}
  \end{equation}
  where $\dual{\cdot}{\cdot}$ on the right hand side denotes   
  the duality between 
  $W^{1-\frac{1}{q'},q'}(\partial\subomega)$ and 
  $W^{-1+\frac{1}{q'},q}(\partial\subomega)$.
  If, in addition, the function $j$ is continuously differentiable on
  $\subomega$ and the partial derivatives have continuous extensions
  to $\overline{\subomega}$, then
  \begin{equation*}
    \int_\subomega  j \cdot \grad\psi  \der{x}
    + 
    \int_\subomega \psi \dive j \der{x}
    = 
    \int_{\partial\subomega} 
    \psi|_{\partial\subomega} 
    \nu \cdot j 
    \der{\sigma_{\subomega}}
    \quad 
    \text{for all $\psi \in W^{1,q'}(\subomega)$,}
  \end{equation*}
  where $\nu$ is the outer unit normal of $\partial\subomega$, and
  $\sigma_\omega$ is the arc--measure on $\partial\subomega$.
\end{prop}
\begin{proof}
  The first statement is a slight generalization, \cf\
  \cite[Lemma~2.4]{kare:current}, of well known results from
  \cite[Ch.~1]{temam}.  The second assertion has been proved in
  \cite[Ch.~5.8]{evans:gariepy92}.
\end{proof}
\begin{thm} 
  \label{thm:curnorm}
  If $(\Varphi,\Fullphi)$ is a solution of van Roosbroeck's system in
  the sense of \defnref{def:vanroos}, and $\subomega\subset\Omega$ is
  an open Lipschitz domain, then there are unique continuous
  functions 
  \begin{math}
    j_{k\nu}\from{
      ]\TNULL,T] 
      \to 
      W^{-1 + \frac {1}{q'},q}(\partial\subomega)
    },
  \end{math}
  $k=1,2$, 
  such that
  \begin{equation}
    \label{eq:gauss}
    \frac{\partial}{\partial t}
    \int_\subomega u_k(t) \der{x}
    =
    \dual{j_{k\nu}(t)}{1} 
    + 
    \int_\subomega 
    \reaction_k(t,\Varphi(t),\Fullphi(t)) \der{x}
    , 
    \quad k =1,2,
  \end{equation}
  where $\dual{\cdot}{\cdot}$ again denotes the duality between
  $W^{1-\frac{1}{q'},q'}(\partial\subomega)$ and
  $W^{-1+\frac{1}{q'},q}(\partial\subomega)$.
\end{thm}
\begin{proof}
  From \eqref{eq:continuity} we obtain for any open Lipschitz domain
  $\subomega\subset\Omega$
  \begin{equation*}
    \int_\subomega 
    u_k'(t) - \dive j_k(t) 
    \der{x}
    =
    \frac{\partial}{\partial t}
    \int_\subomega u_k(t) 
    \der{x}
    -
    \int_\subomega \dive j_k(t)
    \der{x}
    =
    \int_\subomega 
    \reaction_k(t,\Varphi(t),\Fullphi(t))
    \der{x} 
    ,   
  \end{equation*}
  where $j_k$ is defined by \eqref{eq:cur-density}. Using
  \proref{prop:currnor} we find for every $t\in]\TNULL,T]$ a unique
  element $j_{k\nu}(t){\in}W^{-1+\frac{1}{q'},q}(\partial\subomega)$
  such that \eqref{eq:gauss} holds. Moreover, continuity passes over
  from the functions \eqref{eq:jk-regularity} to the mappings
  \begin{math}
    ]\TNULL,T] \ni t 
    \mapsto 
    j_{k\nu}(t) \in W^{-1+\frac{1}{q'},q}(\partial\subomega)
    .
  \end{math}
\end{proof}
If the currents $j_k(t)$ are continuously differentiable on
$\subomega$ and the partial derivatives have continuous extensions to
$\overline{\subomega}$, then by the second part of
\proref{prop:currnor} the formula \eqref{eq:gauss} takes the form
\eqref{eq:balancelaw}.

\section[Numerics]{Numerics}

\thmref{thm:curnorm} is the basis for space discretization of
drift--diffusion equations by means of the finite volume method (FVM).
The FVM was adopted for the numerical solution of van Roosbroeck's
equations by Gajewski, and this approach has been further investigated
in \cite{gajewski:gaertner94,% 
0978.65081,%%Fuhrmann and Langmach
gklnr:detectors,05013697}. %% Eymard/Fuhrmann/Gartner
To discretise the spatial domain one uses a partition into simplex
elements.
Let $\mathcal{E}$ be the set of all edges $e_{il}=x_i-x_l$ of this
triangulation, where $x_1$, $x_2$,\ldots are the vertices.
Moreover, we define the Voronoi cell assigned to a vertex $x_i$ by
\begin{multline*}
  V_i
  \df
  \{
  \text{$x$ in the spatial simulation domain, such that} 
  \\
  \norm{x-x_i} \le \norm{x-x_l}
  \quad
  \text{for all vertices $x_l$ of the triangulation}
  \}
  ,
\end{multline*}
where $\norm{\cdot}$ refers to the norm in the spatial simulation
space $\mathbb{R}^2$.
Now, to get a space discrete version of the current--continuity 
equation, we specify \eqref{eq:gauss} with $\subomega=V_i$, and 
approximate $\dual{j_{k\nu}(t)}{1}$ piecewise by 
\begin{math}
  {j_k}_{il}
  \sigma(\partial V_i \cap \partial V_l),
\end{math}
$\sigma$ being the arc measure on the boundary of $\subomega=V_i$.
The intermediate value ${j_k}_{il}$ can be obtained as follows: The
main hypothesis with respect to the discretization of the currents ---
due to Scharfetter and Gummel \cite{scharfetter:gummel69} --- is that
the electron and hole current density $j_2$ and $j_1$ are constant
along simplex edges.  This assumption allows to calculate ${j_1}_{il}$
and ${j_2}_{il}$ --- the constant values on the edge $e_{il}$ --- in
terms of the node values of the electrostatic potential and the
particle densities, \cf\ \eg\ \cite{gklnr:detectors}.
Thus, one ends up with the following FVM discretization of van
Roosbroeck's system for all interior Voronoi cells $V_i$:
\begin{equation*}
  \begin{split}
    \varepsilon(x_i) 
    \sum_{l \with e_{il}\in\mathcal{E}}
    (\grad\varphi)_{il}
    \sigma(\partial V_k \cap \partial V_l)
    & =
    \left(
      \doping(x_i) + {u_1}(x_i) - {u_2}(x_i)
    \right)
    \abs{V_i},
    \\
    \frac{\partial{u_k}}{\partial t}(x_i)\abs{V_i}
    -{j_k}_{il}\sigma(\partial V_i \cap \partial V_l)
    & =
    \reaction_k(t,\Varphi,\Fullphi_1,\Fullphi_2)(x_i)
    \abs{V_i}
    ,
  \end{split}
\end{equation*}
where $\abs{V_i}$ is the volume of the Voronoi cells $V_i$.
Here we have tested the Poisson equation also with the characteristic
function $1_{V_i}$ of the Voronoi cell $V_i$, and we have applied
Gauss' theorem. In view of \proref{prop:currnor} we assume,
additional to \assuref{assu:dop},
\begin{math}
  \doping
  \from{[\TNULL,\TEINS] \to \widehat{L}^p},
\end{math}
and observe that $\rbvarphi$ can be choosen such that
\begin{math}
  \dual{\rbvarphi}{1_{V_i}}=0
\end{math}
for interior Voronoi cells $V_i$,
\cf\ \remref{r-extend}.
Again, we approximate the right hand side of \eqref{eq:currnor} 
piecewise by
\begin{math}
  (\grad\varphi)_{il} \sigma(\partial V_i \cap \partial V_l),
\end{math}
and we assume --- in consonance with the hypothesis about currents --- 
that the gradient of the electrostatic potential is constant
on the edges of the triangulation, that means
\begin{math}
  (\grad\varphi)_{il}
  =
  (\varphi(x_i)-\varphi(x_l))/\norm{x_i-x_l}%_{\mathbb{R}^2}
  .
\end{math}

Usually, this finite volume discretization of space has been combined
with implicit time discretization, \cf\ \eg\ \cite{gajewski93}. Please
note that the strong differentiability of the electron and hole
density in time is constitutive for this approach. 
% 
% On the other hand,
% if additionaly a Lyapunov function exists, then it provides an
% estimate of the step size in time, \cite{gajewski93}.

\section[Outlook to three spatial dimensions]%
{Outlook  to three spatial dimensions}
\label{sec:outlook}

Much of semiconductor device simulation relies on spatially
two-dimensional models. However, with increasing complexity of
electronic device design spatially three-dimensional simulations
become ever more important, \cf\ \eg\
\cite{gklnr:detectors,gartner:richter:06,gartner:depfet}.  This raises
the question which of the results for the two-dimensional case carry
over to the three-dimensional case. In particular, can one expect that
in three spatial dimensions the divergence of the currents belongs to
a Lebesgue space, and is it possible to establish strong
differentiability of the carrier densities under the rather weak
assumptions about the reaction terms of this paper.

Conditio sine qua non for a modus operandi as in this paper is that in
the three-dimensional case the operators
\begin{equation*}
  -\nabla \cdot \varepsilon \nabla 
  \from{
    \widehat{W}^{1,q}_{\widehatNeumann}
    \to
    \widehat{W}^{-1,q}_{\widehatNeumann}
  }
  \quad
  \text{and} 
  \quad
  -\nabla \cdot \mu_k \nabla 
  \from{
    {W}^{1,q}_{\Neumann}
    \to
    {W}^{-1,q}_{\Neumann}
  }
\end{equation*}
provide isomorphisms for a summability index $q>3$. Unfortunately,
this is not so for arbitrary three-dimensional spatial domains, \cf\
\cite{meyers63}. However, one can proof such a result for certain
classes of three-dimensional material structures and boundary
conditions, \cf\ \cite{pub:1066}, for instance for layered media and
Dirichlet boundary conditions. Dauge proved the result in
\cite{0767.46026} for the Dirichlet Laplacian on a convex polyhedron,
provided the Dirichlet boundary part is separated from its complement
by a finite union of line segments. It would be satisfactory to
combine this conclusion with a heterogeneous material composition.

Under the hypothesis the afore mentioned isomorphisms exist there are
results on quasilinear parabolic systems --- analogous to
\proref{prop:oint} --- \cf\ \cite{rehberg:2005amann} and
\cite{hieber:rehberg:06}, such that one can obtain classical solutions
of the spatially three-dimensional drift--diffusion equations very
much in the same way as here in the two-dimensional case.

\begin{acknowledgement}
  We would like to thank Klaus G\"artner for discussions about 
  % the numerics of 
  van Roosbroeck's system.
\end{acknowledgement}

\providecommand{\bysame}{\leavevmode\hbox to3em{\hrulefill}\thinspace}
\providecommand{\MR}{\relax\ifhmode\unskip\space\fi MR }
\providecommand{\MRhref}[2]{%
  \href{http://www.ams.org/mathscinet-getitem?mr=#1}{#2}
}
\providecommand{\href}[2]{#2}

% \bibliographystyle{amsplain}
% \bibliography{/Home/ima23/kaiser/tex/bib/common}

\end{document}